\documentclass[]{interact}

\usepackage{epstopdf}
\usepackage[caption=false]{subfig}

\usepackage[numbers,sort&compress]{natbib}
\bibpunct[, ]{[}{]}{,}{n}{,}{,}
\renewcommand\bibfont{\fontsize{10}{12}\selectfont}

\theoremstyle{plain}
\newtheorem{theorem}{Theorem}[section]
\newtheorem{lemma}[theorem]{Lemma}
\newtheorem{corollary}[theorem]{Corollary}
\newtheorem{proposition}[theorem]{Proposition}
\newtheorem{assumption}{Assumption}

\theoremstyle{definition}
\newtheorem{definition}[theorem]{Definition}
\newtheorem{example}{Example}

\theoremstyle{remark}
\newtheorem{remark}{Remark}

\newcommand{\Curly}[1]{\left\{ #1 \right\}}
\newcommand{\D}{\mathcal{D}_\alpha}

\newcommand{\Exp}[1]{\mathbb{E}\left[ #1 \right]}
\newcommand{\ExpP}[1]{\mathbb{E}_P \left[ #1 \right]}
\newcommand{\Ip}[2]{\left\langle #1,#2 \right\rangle}

\newcommand{\Paren}[1]{\left( #1 \right)}
\newcommand{\prob}{{P}}

\newcommand{\Sum}[3]{\sum_{#1}^{#2}{#3}}
\newcommand{\Square}[1]{\left[ #1 \right]}

\newcommand{\Tr}[1]{\text{tr}(#1)}

\newcommand{\lowxi}{\xi^l}
\newcommand{\uppxi}{\xi^u}
\newcommand{\realgap}[1]{f^{ij}_{\mathrm{#1}}}
\newcommand{\meangap}[1]{\text{mean}\ f^i_{\mathrm{#1}}}
\newcommand{\medgap}[1]{\text{med}\ f^i_{\mathrm{#1}}}
\newcommand{\vargap}[1]{\text{sd}\ f^i_{\mathrm{#1}}}

\usepackage{amsmath}
\usepackage{amssymb}
\usepackage{amsthm}
\usepackage{mathrsfs} 
\usepackage{mathtools}
\usepackage[caption=false]{subfig}

\begin{document}

\articletype{original article}

\title{Distributionally Robust Expected Residual Minimization for Stochastic Variational Inequality Problems}

\author{
\name{Atsushi Hori\textsuperscript{a}\thanks{Corresponding author: Atsushi Hori. hori@amp.i.kyoto-u.ac.jp}, Yuya Yamakawa\textsuperscript{a}, and Nobuo Yamashita\textsuperscript{a}}
\affil{\textsuperscript{a}Department of Applied Mathematics and Physics, Graduate School of Informatics, Kyoto University, Japan}
}

\maketitle

\begin{abstract}
    The stochastic variational inequality problem (SVIP) is an equilibrium model that includes random variables and has been widely applied in various fields such as economics and engineering.
    Expected residual minimization (ERM) is an established model for obtaining a reasonable solution for the SVIP, and its objective function is an expected value of a suitable merit function for the SVIP.
    However, the ERM is restricted to the case where the distribution is known in advance.
    We extend the ERM to ensure the attainment of robust solutions for the SVIP under the uncertainty distribution (the extended ERM is referred to as distributionally robust expected residual minimization (DRERM), where the worst-case distribution is derived from the set of probability measures in which the expected value and variance take the same sample mean and variance, respectively).
    Under suitable assumptions, we demonstrate that the DRERM can be reformulated as a deterministic convex nonlinear semidefinite programming to avoid numerical integration.
\end{abstract}


\begin{amscode}
90C33; 90C15; 65K15
\end{amscode}

\section{Introduction}

In this study, we consider the following variational inequality problem with a random vector $\xi$: Find $x^*\in S$ such that
\begin{align}
    & \Ip{F(x^*,\xi)}{x-x^*}\geq 0\ \ \forall x\in S,\label{SVIP} \\
    & \quad \text{for $\xi\in\Xi$, almost surely} \nonumber,
\end{align}
where $F\colon\Re^n\times\Xi\rightarrow\Re^n$, and $S\subset\Re^n$ is closed and convex set.
Hereafter, we consider the case where the probability distribution of the random vector $\xi$ may be unknown and provided only partial information, and let $\Xi\subseteq\Re^m$ be a closed convex set referred to as the support of distributions of $\xi$.
System \eqref{SVIP} is known as the stochastic variational inequality problem (SVIP), and is applied in several fields such as economics or engineering to design a market or traffic model, respectively.
In particular, when the set $S$ is given as the nonnegative orthant $\Re^n_+:=\Curly{x\in\Re^n\mid x\geq 0}$, SVIP \eqref{SVIP} can be deduced as the stochastic nonlinear complementarity problem (SNCP): Find $x^*$ such that
\begin{align}\label{SNCP}
	x^*\geq 0,\ F(x^*,\xi)\geq 0,\ \Ip{F(x^*,\xi)}{x^*}=0,
\end{align}
and it has also been studied for a long time.
If the mapping $F$ is linear, SNCP \eqref{SNCP} is referred to as the stochastic linear complementarity problem (SLCP).

In general, there may be no solution that satisfies \eqref{SVIP} or \eqref{SNCP} for almost every $\xi\in\Xi$; thus, the important goal is to find a reasonable solution that minimizes the violation of \eqref{SVIP}.
To obtain such solutions, several models have been proposed such as the expected value (EV) model, expected residual minimization (ERM) model, and distributionally robust model.

The EV model \cite{Gurkan1999} considers the following deterministic variational inequality:
\begin{equation}\label{EVVIP}
    \langle\hat{F}(x^*),x-x^*\rangle\geq 0\quad\forall x\in S,
\end{equation}
where $\hat{F}(x):=\mathbb{E}\Square{F(x,\xi)}$.
Note that an alternative way can also be considered for the expected value of $F$, such as $F(x,\Exp{\xi})$; however, this is not equivalent to the mapping $\hat{F}$ in general.

On the other hand, the ERM was proposed by Chen and Fukushima \cite{Chen2005} for the SLCP.
The primary purpose was to reformulate \eqref{SNCP} as a stochastic optimization problem by using a merit function for the LCP, e.g., the squared Fischer--Burmeister function.
They verified that the ERM tends to output more conservative solutions compared with the EV because the ERM is designed to minimize the mean distance to the solution set of VI for each $\xi\in\Xi$, while the EV only considers the mean $\hat{F}$ of the mapping ${F}(\cdot,\xi)$.

As the natural extension, the ERM for SVIP \eqref{SVIP} can be considered as follows by using a merit function $f(\cdot,\xi)\colon\Re^n\rightarrow\Re_+$ for variational inequalities: 
\begin{align*}
    \begin{array}{rcl}
        \mathrm{(ERM)} & \min & \mathbb{E}\Square{f(x,\xi)} \\
        & \text{s.t.} & x\in S,
    \end{array}
\end{align*}
where the function $f(\cdot,\xi)$ satisfies the following properties for any fixed $\xi\in\Xi$:
\begin{itemize}
    \item[(i)]  $f(x,\xi)\geq 0$ for every $x\in S$;
    \item[(ii)] $x^*\in S$ is a solution of the VIP if and only if $f(x^*,\xi)=0$.
\end{itemize}
To date, several ERM models have been proposed corresponding to each merit function \cite{Luo2009a,Luo2009b,Chen2012}.

However, the ERM has two drawbacks.
First, its distribution of the random vector $\xi$ is assumed to be known in spite of the fact that it may not be observed in various real situations.
Even if one can estimate a distribution from observations, the reliability and robustness of solutions for SVIP \eqref{SVIP} or SNCP \eqref{SNCP} are not guaranteed unless the estimation is sufficiently close to the true distribution, which is referred to as `black swans' in risk theory.
Second, the ERM requires a numerical integration such as the (quasi-)Monte Carlo method to evaluate the expected residual value.
However, the numerical integration is computationally expensive in general; it is advisable to avoid such a sample-based approach.

To tackle these issues, Zhu et al. \cite{Zhu2017} proposed the following conservative approximation model for SNCP \eqref{SNCP}: 
\begin{align}\label{prob:ZhuDRCP}
	\begin{array}{cl}
		\underset{x\in\Re^n}{\min} & \displaystyle\sup_{\prob\in\mathscr{P}}\ \Curly{\ExpP{\Psi(x,\xi)} \mid P(\{F(x,\xi)\geq 0\}\cap\Xi)\geq 1-\varepsilon} \\
		\text{s.t.} & x\geq 0,
	\end{array}
\end{align}
where $0<\varepsilon<1$ is a tolerance parameter, and $\Psi\colon\Re^n\times\Xi\rightarrow\Re$ is a complementarity measure, e.g., $\Psi(x,\xi)=\|x\circ F(x,\xi)\|_2^2$, where $\circ$ denotes the Hadamard product defined by $x\circ y=(x_1y_1,x_2y_2,\ldots,x_ny_n)$ for the vectors $x\in\Re^n$ and $y\in\Re^n$.
Here, $\ExpP{\cdot}$ is the expected value with respect to a distribution function $P(\cdot)\in\mathscr{P}$, where $\mathscr{P}$ is an uncertainty set of the distribution functions supported over $\Xi$ called an ambiguity set.
They considered $\mathscr{P}$ as the following moment ambiguity set:
\begin{align}\label{ambiguity0}
	\mathscr{P}=\left\{
		P\in\mathscr{M}_\Xi\ \middle|\ 
		\mathbb{E}_P[\xi]=\mu_0,\mathbb{E}_P[\xi\xi^\top]=\Sigma_0+\mu_0\mu_0^\top
	\right\},
\end{align}
where $\mathscr{M}_\Xi$ denotes a set of all probability measures supported over $\Xi$, and $\mu_0$ and $\Sigma_0$ respectively denote the (estimated) mean and variance of $\xi$ from observation.
Then they reformulated \eqref{prob:ZhuDRCP} into a nonlinear semidefinite programming problem (NSDP).
In the definition of \eqref{ambiguity0}, however, it is implicitly assumed that an observer knows the exact mean $\mu_0$ and variance $\Sigma_0$.
In the absence of this assumption, the model may not perform properly because observation errors are not considered.
In terms of the distributionally robust optimization (DRO), it is often considered that $\mu_0$ and $\Sigma_0$ cannot be estimated exactly, e.g., the lack of sample data, which motivates us to adopt a more general moment ambiguity set.

In this study, we propose a distributionally robust model of SVIP \eqref{SVIP} under uncertainty of distribution, where the ambiguity set is based on Delage and Ye \cite{Delage2010} (eq. \eqref{eq:moment_set} in Assumption~\ref{as:positive_definiteness}).
Note that our methodology differs from an analysis of the (qualitative or quantitative) statistical robustness \cite{Jiang2022,Guo2021,Kratschmer2012} of a solution obtained from a sample average approximation approach, whose data may contain noise; this is one of the key concepts to study a stochastic model under the uncertainty distribution.
This paper rather focuses on distributional robustness by constructing the ambiguity set with the data.
We propose the following distributionally robust ERM (DRERM) model:
\begin{align*}
    \begin{array}{rcl}
        \mathrm{(DRERM)} & \displaystyle \min & \displaystyle \sup_{{P}\in\mathscr{P}}\ExpP{f(x,\xi)} \\
        & \text{s.t.} & x \in S.
    \end{array}
\end{align*}
This model can be regarded as an extension of the ERM and utilizes some remarkable aspects as stated below:
We illustrate a reformulation of (DRERM) into an NSDP under certain suitable assumptions.
Consequently, it is not required to compute numerical integrals to evaluate the expected value of the stochastic gap functions.

In this paper, we mainly focus on the following regularized gap function \cite{Fukushima1992} as a merit function $f$ in (DRERM):
\begin{align}\label{eq:regularized_gap}
    f(x,\xi)=f_\alpha(x,\xi)\coloneqq\max_{y\in S}\left\{\Ip{F(x,\xi)}{x-y} - \frac{1}{2\alpha}\|y-x\|^2 \right\}.
\end{align}
where $\alpha>0$ is a regularization parameter.
When $S=\Re^n$, the regularized gap function is reduced to $(\alpha/2) \|F(x, \xi)\|^2$.
Therefore, the ERM with $f_\alpha$ is regarded as an extension of the least square problem, and hence it is popular \cite{Agdeppa2010,Luo2009a,Luo2009b,Chen2017}.
Moreover, as we will see in Section \ref{ssec:convexity}, (DRERM) with $f_\alpha$ can be reformulated into a convex NSDP for certain SVIPs.
Note that the NSDP approximation proposed in \cite{Zhu2017} is not convex in general.

The remainder of this paper is organized as follows.
In Section \ref{sec:DRERM_reformulation}, we propose an NSDP model that conservatively approximates (DRERM).
In addition, we show the convexity of the NSDP under certain assumptions.
In Section \ref{sec:numerical_experiments}, we conduct two types of numerical experiments to illustrate the behavior of our reformulation model.
In Section \ref{sec:conclusion}, we conclude this study.

Throughout this paper, we use the following notation.
Let $\Ip{X}{Y}:=\Tr{XY}=\Sum{i,j=1}{m}{X_{ij}Y_{ij}}$ be the matrix inner product of $X\in\mathbb{S}^m$ and $Y\in\mathbb{S}^m$, where $\mathbb{S}^m$ is the set of symmetric matrices included in $\Re^{m\times m}$.
If $X$ and $Y$ are column vectors, $\Ip{X}{Y}$ is the Euclidean inner product.
Let $\mathbb{S}^m_{++}\ (\mathbb{S}^m_+)$ be the set of positive (semi)definite matrices on $\mathbb{S}^m$.

\section{Reformulation and convexity of distributionally robust ERM}\label{sec:DRERM_reformulation}
First, we introduce several approaches to solve (DRERM).
Second, we reformulate (DRERM) into a deterministic NSDP to find its solution efficiently.
Finally, we provide a sufficient condition for the convexity of the NSDP when the mapping $F$ is affine with respect to $x$.

A general technique for solving (DRERM), regardless of the definition of $\mathscr{P}$, is to reformulate it into the following semi-infinite programming and apply the cutting-surface method \cite{Mehrotra2014}:
\begin{align}\label{prob:direct_problem}
	\begin{array}{cl}
		\underset{x,\theta}{\min} & \theta \\
		\text{s.t.} & \theta \geq \ExpP{f_\alpha(x,\xi)}\quad\forall\prob\in\mathscr{P}, \\
								& x\in S.
	\end{array}
\end{align}
Moreover, when $\Xi$ is a finite sample space, i.e., $\Xi\coloneqq\Curly{\xi^1,\xi^2,\ldots,\xi^L}$, problem \eqref{prob:direct_problem} is consequently reduced to the following robust optimization problem because $\mathscr{P}$ can be regarded as a subset of $\Re^L$:
\begin{align}\label{prob:direct_robust}
	\begin{array}{cl}
		\underset{x,\theta}{\min} & \theta \\
		\text{s.t.} & \theta \geq\displaystyle \frac{1}{L}\Sum{k=1}{L}{f_\alpha(x,\xi^k)P_k}\quad\forall P\in\mathscr{P}\subset\{P\in\Re^L_+\mid\textstyle\sum_{k=1}^L P_k=1\}, \\
								& x\in S.
	\end{array}
\end{align}
Thus, nonlinear robust optimization frameworks can be directly applied to \eqref{prob:direct_robust}.
For more details, see \cite{Ben-Tal2000,Ben-Tal2002,Bertsimas2006}.

Another strategy to solve (DRERM) is to consider the duality of the inner supremum part:
\begin{equation}\label{prob:sup_part}
    \sup_{{P}\in\mathscr{P}}\ExpP{f_\alpha(x,\xi)},
\end{equation}
and solve the dual problem.
We adopt this approach and demonstrate that (DRERM) can be reformulated as a deterministic NSDP under certain assumptions. 
For more detailed techniques to deal with general DRO, see \cite{Rahimian2019}.

In the remainder of this study, we assume that $\xi$ is a continuous random variable, and the ambiguity set $\mathscr{P}$ is assumed to be given as the following moment set \cite{Delage2010}, which has been widely applied in existing literature on DRO.
\begin{assumption}[Delage and Ye \cite{Delage2010}]\label{as:positive_definiteness}
	The ambiguity set $\mathscr{P}$ is given by
	\begin{align}\label{eq:moment_set}
		&\mathscr{P}\coloneqq\Curly{
			P\in\mathscr{M}_\Xi\ \middle|\ 
			\begin{array}{l}
				\displaystyle\Paren{\ExpP{\xi}-\mu_0}^\top\Sigma_0^{-1}\Paren{\ExpP{\xi}-\mu_0}\leq\gamma_1 \\
				\displaystyle\ExpP{\Paren{\xi-\mu_0}\Paren{\xi-\mu_0}^\top}\preceq\gamma_2\Sigma_0
			\end{array}
		}
	\end{align}
	where $\gamma_1\geq 0$, $\gamma_2\geq 1$, $\mu_0\in\Xi$, and $\Sigma_0\in\mathbb{S}^m_{++}$.
\end{assumption}

The first condition of \eqref{eq:moment_set}, i.e., $\Paren{\ExpP{\xi}-\mu_0}^\top\Sigma_0^{-1}\Paren{\ExpP{\xi}-\mu_0}\leq\gamma_1$, represents the uncertainty of the true mean $\ExpP{\xi}$ given by an ellipsoid centered on the estimated mean $\mu_0$. In addition, if $\gamma_1=0$, then $\ExpP{\xi}=\mu_0$.
The second condition $\mathbb{E}_P[(\xi-\mu_0)(\xi-\mu_0)^\top]\preceq\gamma_2\Sigma_0$ refers to the uncertainty of the true variance-covariance $\mathbb{E}_P[(\xi-\mu_0)(\xi-\mu_0)^\top]$.
The parameters $\gamma_1$ and $\gamma_2$ determine the strength of the confidence of estimations $\mu_0$ and $\Sigma_0$, respectively; hence, they are referred to as confidence parameters.
A method for determining suitable $\gamma_1$ and $\gamma_2$ from observed samples is introduced in Section 3.4 in~\cite{Delage2010}.
	
\begin{remark}
	When $\gamma_1=0$, $\gamma_2=1$, and the equality holds in the variance-covariance condition in \eqref{eq:moment_set}, the set $\mathscr{P}$ is reduced to \eqref{ambiguity0} considered by Zhu et al {\rm\cite{Zhu2017}}.
\end{remark}

Under Assumption \ref{as:positive_definiteness}, we obtain the following property.

\begin{theorem}\label{thm:sip_reformulation}
	Suppose that Assumption {\rm\ref{as:positive_definiteness}} holds.
	Then $\mathrm{(DRERM)}$ is equivalently reformulated as the following semi-infinite programming with second-order cone constraints:
    \begin{align*}
        \begin{array}{rcl}
            \mathrm{(SIP)} & \underset{(x,y_0,y,Y,z_0)\in\mathcal{V}}{\min} & y_0 + z_0 + \mu_0^\top y + \Ip{\gamma_2\Sigma_0+\mu_0\mu_0^\top}{Y} \\
            & \mathrm{s.t.} & z_0 \geq \sqrt{\gamma_1}\left\|\Sigma_0^{1/2}(y+2Y\mu_0)\right\|,\\
            & & \xi^\top Y \xi + \xi^\top y + y_0 \geq f_\alpha(x,\xi)\quad\forall \xi\in\Xi, \\
            & & x\in S,\ Y\in\mathbb{S}^{m}_+,
        \end{array}
    \end{align*}
    where $\mathcal{V}\coloneqq\Re^n\times\Re\times\Re^m\times\mathbb{S}^m\times\Re$.
\end{theorem}
\proof
From Assumption \ref{as:positive_definiteness} and Lemma 1 of \cite{Delage2010}, for any fixed $x$, the optimal value of \eqref{prob:sup_part}, which is denoted by $\Psi(x;\gamma_1,\gamma_2)$, is equal to that of the following dual problem of \eqref{prob:sup_part}:
\begin{align}\label{prob:dual_sup_part}
		\begin{array}{cl}
			\underset{y_0,y,Y,z_0}{\min} & y_0 + z_0 + \mu_0^\top y + \Ip{\gamma_2\Sigma_0+\mu_0\mu_0^\top}{Y} \\
							\mathrm{s.t.} & z_0 \geq \sqrt{\gamma_1}\left\|\Sigma_0^{1/2}(y+2Y\mu_0)\right\|,\\
							& \xi^\top Y \xi + \xi^\top y + y_0 \geq f_\alpha(x,\xi)\quad\forall \xi\in\Xi, \\
						    & Y\in\mathbb{S}^{m}_+.
		\end{array}
\end{align}
Thus, we obtain the equivalent reformulation of (DRERM) by considering $\min$\ $\{\Psi(x;\gamma_1,\gamma_2)\mid x\in S\}$.
Since optimal values of (SIP) and (DRERM) are equal, the assertion is shown.
\qed

\subsection{Reformulation of SIP into NSDP}\label{ssec:transform_sip}

The goal of this section is to prove that the semi-infinite constraint
\begin{align}\label{ieq:semi-infinite_constraint}
    \xi^\top Y\xi + \xi^\top y + y_0 \geq f_\alpha(x,\xi)\quad\forall\xi\in\Xi,
\end{align}
can be reformulated as a semidefinite constraint by using the duality for the inner maximization of \eqref{eq:regularized_gap}.

In the remainder of this paper, we assume that the closed convex set $S$ is given as a polyhedron:
\[
    S\coloneqq\{x\in\Re^n\mid Ax=b, x\geq 0\},
\]
where $A\in\Re^{l \times n}$ and $b\in\Re^l$.

First, we provide an equivalent form of \eqref{ieq:semi-infinite_constraint} by using the strong duality of the maximization problem in \eqref{eq:regularized_gap}.
\begin{lemma}\label{lem:equivalent_form_of_semi-infinite_constraint}
    The point $(x,y_0,y,Y)$ satisfies \eqref{ieq:semi-infinite_constraint} if and only if there exists $(\lambda,\mu)\in\Re^l\times\Re^n_+$ such that
    \begin{align}\label{ieq:transformed_semi-infinite_constraint}
        \xi^\top Y\xi + \xi^\top y + y_0 \geq \omega_\alpha(x,\lambda,\mu;\xi)\quad\forall\xi\in\Xi.
    \end{align}
    Here,
    \begin{align}\label{eq:omega_function}
        \omega_\alpha(x,\lambda,\mu;\xi)\coloneqq\frac{\alpha}{2}\|F(x,\xi)+A^\top\lambda-\mu\|^2+\Ip{b-Ax}{\lambda}+\Ip{\mu}{x}.
    \end{align}
\end{lemma}
\proof
First, we prove \eqref{ieq:transformed_semi-infinite_constraint} implies \eqref{ieq:semi-infinite_constraint}.
We have the following minimization problem by considering the duality of the maximization problem included in \eqref{eq:regularized_gap}.
\begin{align}\label{prob:RGAPdual_prob}
	\begin{array}{cl}
		\underset{(\lambda,\mu)\in\Re^l\times\Re^n}{\min} & \omega_\alpha(x,\lambda,\mu;\xi) \\
		\text{s.t.} & \mu\in\Re^n_+
	\end{array}
\end{align}
From the weak duality, we have $\omega_\alpha(x,\lambda,\mu;\xi)\geq f_\alpha(x,\xi)$ for each $(x,\xi)\in S\times\Xi$.
Thus, if there exists $(\lambda,\mu)\in\Re^l\times\Re^n_+$ such that $(x,\lambda,\mu,y_0,y,Y)$ satisfies \eqref{ieq:transformed_semi-infinite_constraint}, then the point $(x,y_0,y,Y)$ satisfies \eqref{ieq:semi-infinite_constraint}.

Next, we prove the converse, i.e., \eqref{ieq:semi-infinite_constraint} implies \eqref{ieq:transformed_semi-infinite_constraint}.
The inner maximization in the function $f_\alpha$ is a convex optimization problem whose optimal value is finite for any $x\in S$.
Moreover, owing to the strong duality, there exists $(\lambda,\mu)\in\Re^l\times\Re^n_+$ such that $f_\alpha(x,\xi)=\omega_\alpha(x,\lambda,\mu;\xi)$ for each $(x,\xi)\in S\times\Xi$.
Therefore, if $(x,y_0,y,Y)$ satisfies the condition \eqref{ieq:semi-infinite_constraint}, there exists $(\lambda,\mu)\in\Re^l\times\Re^n_+$ such that \eqref{ieq:transformed_semi-infinite_constraint} holds.
\qed
\bigskip

Now, we make assumptions on the mapping $F$ and the support $\Xi$ in SVIP \eqref{SVIP}.
Similar assumptions on $F$ and $\Xi$ have already been considered by Zhu et al. \cite{Zhu2017} for SNCP \eqref{SNCP}\footnote{Only when the complementarity measure is evaluated by $\|x\circ F(x,\xi)\|_\infty$, the mapping $F$ of \eqref{prob:ZhuDRCP} is allowed up to second-order with respect to $\xi$.}.
For certain examples that satisfy the following assumptions on SVIP \eqref{SVIP}, see \cite{Agdeppa2010}.
\begin{assumption}\label{asmp:FandS}\mbox{\\}
    \begin{enumerate}
        \item [$\mathrm{(i)}$] The $i$-th element of the mapping $F$ is affine with respect to $\xi$:
    	\begin{align*}
    		F_i(x,\xi)\coloneqq \Paren{c^i(x)}^\top\xi + c^i_0(x),\ i=1,2,\ldots,n.
    	\end{align*}
        \item [$\mathrm{(ii)}$] The support $\Xi$ is given as
        \begin{align}\label{eq:Xi}
        	\Xi\coloneqq\Curly{\xi\in\Re^m\mid g_i(\xi)\leq 0, i=1,2,\ldots,p}.
        \end{align}
        Here, $g_i\colon\Re^m\to\Re$ is defined by
        \begin{align}\label{eq:g_i_quad_form}
            g_i(\xi)\coloneqq\xi^\top\tilde{A}_i\xi+2{\tilde{b}_i}^\top\xi+\tilde{c}_i,\quad i=1,2,\dots,p,
        \end{align}
        where $\tilde{A}_i\in\mathbb{S}^m$, $\tilde{b}_i\in\Re^m$, and $\tilde{c}_i\in\Re$.
    \end{enumerate}
\end{assumption}

As preliminaries, let us introduce the S-procedure and its special case.

\begin{lemma}[S-procedure Derinkuyu and Pınar \cite{Derinkuyu2006}]\label{def:S-procedure}
	Let $\Xi$ be given as \eqref{eq:Xi} and
	\begin{align}\label{g0func}
	    g_0(\xi)\coloneqq\xi^\top \bar{A}_0\xi+2\xi^\top \bar{b}_0+\bar{c}_0,
	\end{align}
	where $\bar{A}_0\in\mathbb{S}^m$, $\bar{b}_0\in\Re^m$, and $\bar{c}_0\in\Re$.
	Assume that there exists $s\in\Re^p_+$ such that 
	\begin{align}\label{ieq:S-procedure}
		g_0(\xi)+\Sum{i=1}{p}{s_i g_i(\xi)}\geq 0\quad \forall\xi\in\Re^m.
	\end{align}
	Then, $g_0(\xi)\geq 0$ for all $\xi\in\Xi$.
\end{lemma}

The following lemma indicates that the converse also holds when $p=1$ in Lemma~\ref{def:S-procedure}.
\begin{lemma}[P\'{o}lik and Terlaky \cite{Polik2007}]\label{lem:S-lemma}
	Suppose that $\Xi$ is given by \eqref{eq:Xi} with $p=1$ and let $g_0(\xi)$ be defined as \eqref{g0func}.
	Assume that there exists $\hat{\xi}_0$ such that $g_1(\hat{\xi}_0)<0$.
	Then, the statements {\rm(i)} and {\rm(ii)} are equivalent:
	\begin{description}
		\item[$\mathrm{(i)}$] For all $\xi\in\Re^m$, $g_1(\xi)\leq 0$ implies $g_0(\xi)\geq 0${\rm;}
		\item[$\mathrm{(ii)}$] there exists some nonnegative number $s\geq 0$ such that
			\[
				g_0(\xi)+s g_1(\xi)\geq 0\quad\forall\xi\in\Re^m.
			\]
	\end{description}
\end{lemma}


We further introduce an equivalence between nonnegative quadratic functions on $\Re^m$ and semidefiniteness.
\begin{lemma}[Proposition 2 in Sturm and Zhang \cite{Sturm2003}]\label{lem:NSD_constraint}
	Let $\tilde{A}\in\mathbb{S}^m,\ \tilde{b}\in\Re^m,$ and $\tilde{c}\in\Re$ be given.
	Then, the following two conditions {\rm(i)} and {\rm(ii)} are equivalent:
    \begin{align*}
        \begin{array}{cl}
            {\mathrm{(i)}} & \left[1,\ \xi^\top\right]
            \left[
                \begin{array}{cc}
                    \tilde{c} & \tilde{b}^\top \\
                    \tilde{b} & \tilde{A}
                \end{array}
            \right]
            \left[
                \begin{array}{c}
                    1 \\
                    \xi
                \end{array}
            \right]
            \geq 0\quad \forall\xi\in\Re^m;\\
            {\mathrm{(ii)}} & 
            \left[
                \begin{array}{cc}
                    \tilde{c} & \tilde{b}^\top \\
                    \tilde{b} & \tilde{A}
                \end{array}
            \right]
            \succeq O.
        \end{array}
    \end{align*}
\end{lemma}

Zhu et al. \cite{Zhu2017} proposed a certain NSDP that conservatively approximates DRO \eqref{prob:ZhuDRCP}, where the conservative approximation denotes that the optimal value of the NSDP is not less than that of DRO \eqref{prob:ZhuDRCP}.
In this paper, we also provide the following conservative approximation of (DRERM) based on their technique.
\begin{align*}
	\begin{array}{rcl}
		\mathrm{(NSDP)} & \underset{(w,z_0,s)\in\mathcal{W}\times\Re\times\Re^p}{\min} & z_0 + y_0 + \mu_0^\top y + \Ip{\gamma_2\Sigma_0+\mu_0\mu_0^\top}{Y} \\
		& \text{s.t.} & z_0 \geq \sqrt{\gamma_1}\left\|\Sigma_0^{1/2}(y+2Y\mu_0)\right\|, \\
							&	& \displaystyle\mathcal{D}_\alpha(w)+\Sum{i=1}{p}{s_i\tilde{\mathcal{A}}_i}\succeq O, \\
							&	& x\in S,\ \mu\in\Re^n_+,\ Y\in\mathbb{S}^{m}_+,\ s\in\Re^p_+,
    \end{array}
\end{align*}
    where $w\coloneqq(x,\lambda,\mu,y_0,y,Y)\in\mathcal{W}\coloneqq\Re^n\times\Re^l\times\Re^n\times\Re\times\Re^m\times\mathbb{S}^m$, and $\D\colon\mathcal{W}\to\mathbb{S}^{m+1}$ is a symmetric-matrix-valued function defined as follows:
\begin{align}
	\D(w)\coloneqq
		\Square{
        \begin{array}{cc}
            y_0 & 1/2 y^\top \\
            1/2 y & Y
        \end{array}
			} - \Curly{G(x,\lambda,\mu) + \frac{\alpha}{2}\Sum{i=1}{n}{H^i(x,\lambda,\mu)}}, \label{eq:matD}
\end{align}
where
\begin{align*}
    G(x,\lambda,\mu)&\coloneqq
    \Square{
        \begin{array}{cc}
            \Ip{b-Ax}{\lambda} + \Ip{\mu}{x} & 0^\top \\
            0 & O_{m \times m}
        \end{array}
    },
		\\
    H^i(x,\lambda,\mu)&\coloneqq
    \Square{
        \begin{array}{cc}
            p^i_0(x,\lambda,\mu)^2 &
            p^i_0(x,\lambda,\mu)c^i(x)^\top \\
            p^i_0(x,\lambda,\mu)c^i(x) &
            c^i(x)c^i(x)^\top
        \end{array}
    },\ i=1,2,\ldots,n,\\
    p^i_0(x,\lambda,\mu)&\coloneqq c^i_0(x)+\Sum{j=1}{l}{a^{ji}\lambda_j-\mu_i},\ i=1,2,\dots,n,
\end{align*}
and $\tilde{\mathcal{A}}_i$ is defined as
\[
    \tilde{\mathcal{A}}_i\coloneqq\Square{\begin{array}{cc}
												\tilde{c}_i & \tilde{b}^\top_i \\
												\tilde{b}_i	& \tilde{A}_i
										\end{array}}\quad i=1,2,\dots,p.
\]

Next, we provide several definitions and lemmas to prove that (NSDP) gives a conservative approximation of (DRERM).
Now, we define
\begin{align*}
    \tilde{A}_0 \coloneqq Y-\frac{\alpha}{2}\sum_{i=1}^n c^i(x)c^i(x)^\top,\ \tilde{b}_0 \coloneqq\frac{y}{2}-\frac{\alpha}{2}\sum_{i=1}^n p^i_0(x,\lambda,\mu)c^i(x),\\
    \tilde{c}_0 \coloneqq y_0-\langle b-Ax,\lambda\rangle-\langle \mu,x\rangle-\frac{\alpha}{2}\sum_{i=1}^n p^i_0(x,\lambda,\mu)^2,
\end{align*}
and
\begin{align}\label{eq:g_0function}
    h(\xi)\coloneqq \xi^\top\tilde{A}_0\xi+2\xi^\top\tilde{b}_0+\tilde{c}_0.
\end{align}

Under Assumption \ref{asmp:FandS}--(i), \eqref{eq:omega_function} is written as
\begin{align}\label{eq:omega_specific_form}
    \omega_\alpha(x,\lambda,\mu;\xi)=\Square{1,\ \xi^\top}\Curly{G(x,\lambda,\mu) + \frac{\alpha}{2}\Sum{i=1}{n}{H^i(x,\lambda,\mu)}}\Square{\begin{array}{c}1\\\xi\end{array}}.
\end{align}
Through the straightforward calculation, we obtain the following equalities:
\begin{align}\label{D_alpha=omega}
    & \Square{1,\ \xi^\top}\D(w)\Square{\begin{array}{c}1\\\xi\end{array}}=\xi^\top Y\xi+\xi^\top y+y_0 - \omega_\alpha(x,\lambda,\mu;\xi)=h(\xi).
\end{align}

\begin{lemma}\label{lem:nsd_equiv_expression}
    The nonlinear semidefinite constraint included in $\mathrm{(NSDP)}$, i.e.,
    \begin{align}\label{semidefinite_constraint_NSDP}
        \D(w)+\sum_{i=1}^p s_i\tilde{\mathcal{A}}_i\succeq O,
    \end{align}
    is equivalent to
    \begin{align}\label{semiinfinite_nonnegative}
        h(\xi)+\sum_{i=1}^p s_i g_i(\xi)\geq 0\quad\forall\xi\in\Re^m.
    \end{align}
\end{lemma}
\proof
    By Lemma~\ref{lem:NSD_constraint}, \eqref{semidefinite_constraint_NSDP} is equivalent to
    \begin{align}\label{SDP_constr_eq}
       \Square{1,\ \xi^\top}\left(\D(w) + \sum_{i=1}^p s_i \tilde{\mathcal{A}}_i\right)\Square{\begin{array}{c}1\\\xi\end{array}}
    	\geq 0\quad \forall\xi\in\Re^m.
    \end{align}
    Since $[1,\ \xi^\top]\D(w)\left[\begin{array}{c}1\\\xi\end{array}\right]=h(\xi)$ from \eqref{D_alpha=omega} and $[1,\ \xi^\top]\tilde{\mathcal{A}}_i\left[\begin{array}{c}1\\\xi\end{array}\right]=g_i(\xi)$, \eqref{SDP_constr_eq} can be equivalently represented as \eqref{semiinfinite_nonnegative}.
\qed

\bigskip

The next lemma provides a sufficient condition for semi-infinite constraint \eqref{ieq:semi-infinite_constraint}.

\begin{lemma}\label{lem:semi-defi_and_semi-inf}
    Suppose that Assumption {\rm \ref{asmp:FandS}} holds.
    Whenever $p\geq 1$, if there exists $(w,s)\in\mathcal{W}\times\Re^p$ such that $\mu\in\Re^n_+$, $s\in\Re^p_+$, and \eqref{semidefinite_constraint_NSDP}, i.e.,
    \begin{align*}
        \D(w)+\sum_{i=1}^p s_i\tilde{\mathcal{A}}_i\succeq O,
    \end{align*}
    then the subvector $(x,y_0,y,Y)$ satisfies the semi-infinite constraint \eqref{ieq:semi-infinite_constraint}, i.e.,
    \begin{align*}
        \xi^\top Y\xi+\xi^\top y+y_0\geq f_\alpha(x,\xi)\quad\forall\xi\in\Xi.
    \end{align*}
    Furthermore, when $p=1$ and the assumption of Lemma {\rm \ref{lem:S-lemma}} holds, the converse is also true, i.e., if $(x,y_0,y,Y)$ satisfies \eqref{ieq:semi-infinite_constraint}, then there exists $(\lambda,\mu,s)\in\Re^l\times\Re^n_+\times\Re^p_+$ such that \eqref{semidefinite_constraint_NSDP} satisfies.
\end{lemma}
\proof
    First, we prove the general case where $p\geq 1$.
    Assume that there exist $w\in\mathcal{W}$ and $s\in\Re^p_+$ such that semidefinite constraint \eqref{semidefinite_constraint_NSDP} holds.
    Then, by Lemma~\ref{lem:nsd_equiv_expression}, we have \eqref{semiinfinite_nonnegative}, i.e.,
    \[
        h(\xi)+\sum_{i=1}^p s_i g_i(\xi)\geq 0\quad\forall\xi\in\Re^m.
    \]
    By regarding $h(\xi)$ as $g_0(\xi)$ in Lemma~\ref{def:S-procedure}, \eqref{semiinfinite_nonnegative} implies $h(\xi)\geq 0$ for all $\xi\in\Xi$, and it then follows from \eqref{D_alpha=omega} that \eqref{ieq:transformed_semi-infinite_constraint} holds, i.e.,
    \[
        \xi^\top Y\xi+\xi^\top y+y_0\geq\omega_\alpha(x,\lambda,\mu;\xi)\quad\forall\xi\in\Xi.
    \]
    Finally, Lemma~\ref{lem:equivalent_form_of_semi-infinite_constraint} states that $w$ satisfies \eqref{ieq:transformed_semi-infinite_constraint} if and only if its subvector $(x,y_0,y,Y)$ satisfies semi-infinite constraint \eqref{ieq:semi-infinite_constraint}.
    The first part of the proof is completed.
    
    Next, we prove that the converse when $p=1$ and the assumption of Lemma~\ref{lem:S-lemma} holds.
    Suppose that $(x,y_0,y,Y)$ satisfies \eqref{ieq:semi-infinite_constraint}.
    By Lemma~\ref{lem:equivalent_form_of_semi-infinite_constraint}, \eqref{ieq:semi-infinite_constraint} holds if and only if there exists $(\lambda,\mu)\in\Re^l\times\Re^n_+$ such that \eqref{ieq:transformed_semi-infinite_constraint} holds.
    Note that under Assumption \ref{asmp:FandS}--(i), the function $\omega_\alpha(x,\lambda,\mu;\xi)$ is given as \eqref{eq:omega_specific_form}.
    Then, $\xi^\top Y\xi+\xi^\top y+y_0-\omega_\alpha(x,\lambda,\mu;\xi)\geq 0$ and \eqref{D_alpha=omega} yield $h(\xi)\geq 0$.
    Note that $h(\xi)\geq 0$ for all $\xi\in\Xi$ if and only if for all $\xi\in\Re^m$, $g_1(\xi)\leq 0$ implies $h(\xi)\geq 0$.
    It then follows from Assumption \ref{asmp:FandS}--(ii) with $p=1$ and Lemma~\ref{lem:S-lemma} that there exists $s\in\Re_+$ such that $h(\xi)+sg_1(\xi)\geq 0$ for all $\xi\in\Re^m$.
    By Lemma~\ref{lem:nsd_equiv_expression}, this condition is equivalent to semidefinite constraint \eqref{semidefinite_constraint_NSDP} in (NSDP).
    Thus, we have proved the converse.
\qed

\bigskip
The following result shows the feasibility between constraints of (SIP) and (NSDP).

\begin{proposition}\label{prop:feasibility_sip_nsdp}
    Suppose that Assumption~{\rm \ref{asmp:FandS}} holds.
    Whenever $p\geq 1$, if $(w,z_0,s)\in\mathcal{W}\times\Re\times\Re^p$ is feasible to {\rm (NSDP)}, then its subvector $(x,y_0,y,Y,z_0)\in\mathcal{V}$ is also feasible to {\rm (SIP)}.
    Moreover, when $p=1$ and the assumption of Lemma {\rm \ref{lem:S-lemma}} holds, if $(x,y_0,y,Y,z_0)\in\mathcal{V}$ is feasible to {\rm (SIP)}, then there exists $(\lambda,\mu,s)\in\Re^l\times\Re^n_+\times\Re^p_+$ such that $(w,z_0,s)\in\mathcal{W}\times\Re\times\Re_+$ is also feasible to {\rm (NSDP)}.
\end{proposition}
\proof
    Note that all constraints in (NSDP) except \eqref{semidefinite_constraint_NSDP} coincide with those in (SIP) excluding semi-infinite constraint \eqref{ieq:semi-infinite_constraint}.
    This statement and Lemma~\ref{lem:semi-defi_and_semi-inf} ensure that if $p\geq 1$, and $(w,s,z_0)$ is the feasible solution of (NSDP), then its subvector $(x,y_0,y,Y,z_0)$ is the feasible solution to (SIP).
    Thus, we showed the general case where $p\geq 1$.
    
    Suppose that $p=1$, and $(x,y_0,y,Y,z_0)\in\mathcal{V}$ is a feasible solution to (SIP).
    As mentioned above, $(w,z_0,s)$ satisfies the constraints of (NSDP) except \eqref{semidefinite_constraint_NSDP}.
    Moreover, Lemma~\ref{lem:semi-defi_and_semi-inf} guarantees that there exists $(\lambda,\mu,s)\in\Re^l\times\Re^n_+\times\Re^p_+$ such that $(w,z_0,s)$ satisfies constraint \eqref{semidefinite_constraint_NSDP}.
    We have completed the proof.
\qed

\bigskip

By using the above lemmas, we show one of the main results.

\begin{theorem}\label{thm:sufficient_condition_for_sip}
	Suppose that Assumptions {\rm\ref{asmp:FandS}} holds.
	Then, {\rm (SIP)} can be conservatively approximated as {\rm (NSDP)}.
\end{theorem}

\proof
Suppose that $(w,z_0,s)\in\mathcal{W}\times\Re\times\Re^p_+$ is a feasible point of (NSDP).
It then follows from Proposition~\ref{prop:feasibility_sip_nsdp} that the subvector $(x,y_0,y,Y,z_0)$ satisfies the constraints of (SIP).
From the above facts, the optimal value of (NSDP) can never be less than that of (SIP).
Therefore, (NSDP) is a conservative approximation of (SIP).
The proof is completed.
\qed

\bigskip 

Here, we provide some examples of $\Xi$ that can be expressed as the intersection of nonnegative quadratic functions.

\begin{example}[Box set]\label{ex:rectangle}
Consider $\Xi$ given by the following box set:
\[
		\Xi\coloneqq\Curly{\xi\in\Re^m\mid \lowxi_i\leq \xi_i\leq \uppxi_i,\ i=1,2,\ldots,m}.
\]
By using a quadratic function, $\lowxi_i\leq \xi\leq \uppxi_i$ can be rewritten as follows:
\[
		g_i(\xi) = \xi_i(\uppxi_i+\lowxi_i)-\uppxi_i\lowxi_i-\xi_i^2 = \left[1,\ \xi^\top\right]T_i\Square{
				\begin{array}{c}
						1 \\
						\xi
				\end{array}
		}\geq 0,
\]
where
\[
		T_i\coloneqq\Square{
				\begin{array}{cc}
						-\uppxi_i\lowxi_i & -\frac{1}{2}\xi_i(\uppxi_i+\lowxi_i)(e^i)^\top \\
						-\frac{1}{2}\xi_i(\uppxi_i+\lowxi_i)e^i & -\tilde{I}_i
				\end{array}
		}.
\]
Here, $e^i\in\Re^m$ is the $i$-th column vector of the identity matrix, and $\tilde{I}_i\in\Re^{m\times m}$ is a matrix whose elements are all zero except the $(i,i)$ entry which is {\rm1}.

This example corresponds to the case where $\tilde{A}_i=\tilde{I}_i,\tilde{b}_i=\frac{1}{2}\xi_i(\uppxi_i+\lowxi_i)e^i$, and $\tilde{c}_i=\uppxi_i\lowxi_i$ in $\mathrm{(NSDP)}$.
\qed
\end{example}

\begin{example}[Ellipsoids]\label{ex:ellipsoids}
	Consider $\Xi$ given by the following ellipsoids:
	\begin{align}\label{eq:ex:ellipsoids}
		\Xi\coloneqq\Curly{\xi\in\Re^m\mid (\xi-\hat{\xi}^i)^\top P_i^{-1}(\xi-\hat{\xi}^i)\leq 1,i=1,2,\ldots,p},
	\end{align}
	where the vector $\hat{\xi}^i\in\Re^m$ is the center of the $i$-th ellipsoid, and the matrix $P_i$ is supposed to be positive definite.
	This example corresponds to the case where $\tilde{A}_i=P^{-1},\tilde{b}_i=-P_i^{-1}\hat{\xi}^i$, and $\tilde{c}_i=(\hat{\xi}^i)^\top P_i^{-1}\hat{\xi}^i - 1$ in $\mathrm{(NSDP)}$.
	\qed
\end{example}

Next, we illustrate the special case of Theorem~\ref{thm:sufficient_condition_for_sip}, which ensures that a solution of (NSDP) solves (DRERM).

\begin{corollary}\label{cor:single-ellipsoid}
	Suppose that $p=1$ in \eqref{eq:ex:ellipsoids}, and that the assumption of Lemma~\ref{lem:S-lemma} holds.
	Then, if $(w,z_0,s)\in\mathcal{W}\times\Re\times\Re$ is a global optimum of $\mathrm{(NSDP)}$, then $(x,y_0,y,Y,z_0)$ and $x$ are also global optima to $\mathrm{(SIP)}$ and $\mathrm{(DRERM)}$, respectively.
	In addition, the optimal value of $\mathrm{(NSDP)}$ is equal to those of $\mathrm{(SIP)}$ and $\mathrm{(DRERM)}$.
\end{corollary}
\proof
Let $(w,z_0,s)$ be a global optimum to (NSDP).
Assume that its subvector $(x,y_0,y,Y,z_0)$ is not a global optimum of (SIP). 
Note that $(x,y_0,y,Y,z_0)$ is the feasible solution to (SIP) from Proposition~\ref{prop:feasibility_sip_nsdp}.  
By the assumption, there exists a feasible solution $(x',y'_0,y',Y',z'_0)$ in (SIP) such that
\begin{align}\label{ieq:sip_objective}
    z'_0+y'_0+\mu_0^\top y' + \Ip{\gamma_2\Sigma_0+\mu_0\mu_0^\top}{Y'}<z_0+y_0+\mu_0^\top y+ \Ip{\gamma_2\Sigma_0+\mu_0\mu_0^\top}{Y}.
\end{align}
Proposition~\ref{prop:feasibility_sip_nsdp} guarantees that if the solution $(x',y'_0,y',Y',z'_0)$ is the feasible point to (SIP), then there exists $(\lambda',\mu',s')\in\Re^l\times\Re^n_+\times\Re_+$ such that $(w',z'_0,s')\in\mathcal{W}\times\Re\times\Re_+$ is the feasible solution to (NSDP), where $w'\coloneqq(x',\lambda',\mu',y'_0,y',Y')\in\mathcal{W}$.
Because the objective functions of (SIP) and (NSDP) coincide, the solution $(w',z'_0,s')$ of (NSDP) also satisfies the inequality \eqref{ieq:sip_objective}.
Hence, it contradicts that $(w,z_0,s)$ is a global optimum to (NSDP).
We have that $(x,y_0,y,Y,z_0)$, which is the subvector of the global optimum $(w,z_0,s)$ of (NSDP), is also the global optimum in (SIP), and optimal values are equal because their objective functions coincide.
Moreover, since (DRERM) is equivalent to (SIP) from Theorem~\ref{thm:sip_reformulation}, $x$ is also a global optimum to (DRERM), and their optimal values are equal.
\qed

\bigskip

In addition, when $\Xi=\Re^m$, then we can show an equivalence between (SIP) (or (DRERM)) and the following NSDP:
\begin{align*}
	\begin{array}{rcl}
		\mathrm{(NSDP')} & \underset{(w,z_0)\in\mathcal{W}\times\Re}{\min} & z_0 + y_0 + \mu_0^\top y + \Ip{\gamma_2\Sigma_0+\mu_0\mu_0^\top}{Y} \\
		& \text{s.t.} & z_0 \geq \sqrt{\gamma_1}\left\|\Sigma_0^{1/2}(y+2Y\mu_0)\right\|, \\
							&	& \displaystyle\mathcal{D}_\alpha(w)\succeq O, \\
							&	& x\in S,\ \mu\in\Re^n_+.
	\end{array}
\end{align*}
To show this property, we prepare a lemma below.

\begin{lemma}\label{lem:Xi=Rm_SDPeqivform}
Let $(x,y_0,y,Y)$ be given.
Then, the following two statements are equivalent:
\begin{description}
    \item [$\mathrm{(i)}$] There exists $(\lambda,\mu)\in\Re^l\times\Re^n_+$ such that $\D(w)\succeq O$;
    \item [$\mathrm{(ii)}$] $Y\in\mathbb{S}^m_+$ and 
    \begin{align}\label{ieq:semi-infinite_constr_Xi=Rm}
        \xi^\top Y\xi + \xi^\top y+ y_0 \geq f_\alpha(x,\xi)\quad\forall\xi\in\Xi=\Re^m.
    \end{align}
\end{description}
\end{lemma}
\proof
    First, we show that (i) implies (ii).
    By Lemma~\ref{lem:NSD_constraint} and the first equality of \eqref{D_alpha=omega}, $\D(w)\succeq O$ if and only if
    \begin{align}
        \xi^\top Y\xi +\xi^\top y+y_0\geq\omega_\alpha(x,\lambda,\mu;\xi)\quad\forall\xi\in\Re^m.\label{ieq:semi-infinite_constraint_Xi=Rm}
    \end{align}
    As we mentioned in Lemma~\ref{lem:equivalent_form_of_semi-infinite_constraint}, $\omega_\alpha(x,\lambda,\mu;\xi)$ is the dual function of the maximization problem in $f_\alpha$.
    This implies that for any $x\in S$ and $\xi\in\Re^m$, $\omega_\alpha(x,\lambda,\mu;\xi)\geq f_\alpha(x,\xi)\geq$~$0$.
    Then, \eqref{ieq:semi-infinite_constraint_Xi=Rm} implies $\xi^\top Y\xi+\xi^\top y+y_0\geq 0$, and by Lemma~\ref{lem:NSD_constraint}, we have
    \[
        \Square{
            \begin{array}{cc}
                y_0 & 1/2y^\top \\
                1/2 y & Y
            \end{array}
        } \succeq O.
    \]
    By the Schur complement, this ensures the positive semidefiniteness of $Y$.
    Furthermore, \eqref{ieq:semi-infinite_constraint_Xi=Rm} implies \eqref{ieq:semi-infinite_constr_Xi=Rm} by Lemma~\ref{lem:equivalent_form_of_semi-infinite_constraint}.
    We have proved the former part of the proof.
    
    Next, we prove that (ii) implies (i).
    Suppose that $Y\in\mathbb{S}^m_+$ and \eqref{ieq:semi-infinite_constr_Xi=Rm} holds.
    Then, by Lemma~\ref{lem:equivalent_form_of_semi-infinite_constraint}, there exists $(\lambda,\mu)\in\Re^l\times\Re^n_+$ such that \eqref{ieq:semi-infinite_constraint_Xi=Rm} holds, and it is immediately observed that $\D(w)\succeq O$.
    Hence, the proof is completed.
\qed

\bigskip
We obtain the relation regarding the feasibility between (SIP) and $\mathrm{(NSDP')}$ by using Lemma~\ref{lem:Xi=Rm_SDPeqivform}.

\begin{proposition}\label{prop:Xi=Rm_feasibility}
    The point $(x,y_0,y,Y,z_0)\in\mathcal{V}$ is a feasible solution to {\rm (SIP)} if and only if there exists $(\lambda,\mu)\in\Re^l\times\Re^n_+$ such that $(w,z_0)\in\mathcal{W}\times\Re$ is also a feasible solution to $\mathrm{(NSDP')}$.
\end{proposition}
\proof
Similar to the proof of Proposition~\ref{prop:feasibility_sip_nsdp}, all the constraints in $\mathrm{(NSDP')}$ except the semidefinite constraint $\D(w)\succeq O$ coincide with those in (SIP) excluding semi-infinite constraint \eqref{ieq:semi-infinite_constraint}.
This statement and Lemma~\ref{lem:Xi=Rm_SDPeqivform} ensure that for given $(x,y_0,y,Y)$, the point $(x,y_0,y,Y,z_0)$ is feasible to (SIP) if and only if there exists $(\lambda,\mu)\in\Re^l\times\Re^n_+$ such that $(w,z_0)$ is feasible to $\mathrm{(NSDP')}$.
\qed
\bigskip

Finally, the optimality between $\mathrm{(NSDP')}$ and (SIP) is obtained as follows.

\begin{theorem}\label{thm:Xi=Rm}
	Suppose that Assumption {\rm \ref{asmp:FandS}}--{\rm (i)} holds, and that $\Xi=\Re^m$.
	If $(w,z_0)$ is a global optimum to $\mathrm{(NSDP')}$, then its subvector $(x,y_0,y,Y,z_0)$ and $x$ are also global optima for {\rm(SIP)} and {\rm(DRERM)}, respectively.
\end{theorem}
\proof
Let $(w,z_0)$ be a global optimum to $\mathrm{(NSDP')}$.
Assume that its subvector $(x,y_0,y,Y,z_0)$ is not a global optimum of (SIP).
Note that $(x,y_0,y,Y,z_0)$ is a feasible solution to (SIP) by Proposition~\ref{prop:Xi=Rm_feasibility}.
Since $(x,y_0,y,Y,z_0)$ is not a global optimum of (SIP), there exists a feasible solution $(x',y'_0,y',Y',z'_0)$ such that
\begin{align}\label{objfunc_inequality_again}
    z'_0+y'_0+\mu_0^\top y' + \Ip{\gamma_2\Sigma_0+\mu_0\mu_0^\top}{Y'}<z_0+y_0+\mu_0^\top y+ \Ip{\gamma_2\Sigma_0+\mu_0\mu_0^\top}{Y}.
\end{align}
Moreover, by Proposition~\ref{prop:Xi=Rm_feasibility}, the feasible solution $(x',y'_0,y',Y',z'_0)$ of (SIP) is also feasible to $\mathrm{(NSDP')}$ for some $(\lambda',\mu')\in\Re^l\times\Re^n_+$.
This statement and \eqref{objfunc_inequality_again} contradict each other; thus, $(x,y_0,y,Y,z_0)$ is a global optimum of (SIP).

Since (DRERM) is equivalent to (SIP), $x$, which is the subvector of the global optimum $(x,y_0,y,Y,z_0)$ of (SIP), is also a global optimum to (DRERM).
Thus, the optimal value of $\mathrm{(NSDP')}$ coincides with those of (SIP) and (DRERM), respectively.
\qed

\begin{remark}\label{remark:optimality_equivalent}
    Zhu et al. \cite{Zhu2017} have only shown a conservative NSDP approximation for problem \eqref{prob:ZhuDRCP}; that is, the subvector $x\geq 0$ of a global optimal solution of the conservative approximated NSDP may not globally solve \eqref{prob:ZhuDRCP} in general.
    However, as Theorem \ref{thm:Xi=Rm} and Corollary \ref{cor:single-ellipsoid} state, if $\Xi$ is $\Re^m$ or a single ellipsoid, the variable $x\in S$ of a global optimal point obtained from (NSDP) or (NSDP$'$) solves (DRERM).
\end{remark}

\subsection{Convexity of NSDP}\label{ssec:convexity}

First, the sufficient condition is presented under which (NSDP) and $\mathrm{(NSDP')}$ are convex.

\begin{assumption}\label{assmp:affine_F}
    The mapping $F$ is affine with respect to $x$, i.e., 
    \[
    	F(x,\xi):=M(\xi)x+q(\xi),
    \]
    where $M\colon\Xi\rightarrow\Re^{n\times n}$ and $q\colon\Xi\rightarrow\Re^n$.
       Here, the $(i,j)$-entry of $M(\xi)$ is denoted by $\Paren{M(\xi)}_{ij}:=\Paren{m^{ij}}^\top\xi+m^{ij}_0$, and the $i$-th element of $q(\xi)$ is $\Paren{q(\xi)}_i:=\Paren{q^i}^\top\xi+q^i_0$, where $m^{ij},q^i\in\Re^m$, $m^{ij}_0,q^i_0\in\Re$.
       Hence, $c^i(x)$ and $c^i_0(x)$ defined in Assumption {\rm\ref{asmp:FandS}} can be rewritten as follows:
    \begin{align*}
    	c^i(x):=q^i + \bar{M}_ix \in\Re^m,\ i=1,2,\dots,n,\\
    	c^i_0(x):=q^i_0 + (\bar{m}^i_0)^\top x \in\Re,\ i=1,2,\dots,n,
    \end{align*}
    where
    \begin{gather*}
        \bar{M}_i:=[m^{i,1},m^{i,2},\dots,m^{i,n}]\in\Re^{m\times n},\quad i=1,2,\dots,n,\\
        \bar{m}^i_0:=[m^{i,1}_0,m^{i,2}_0,\dots,m^{i,n}_0]^\top\in\Re^n,\quad i=1,2,\dots,n.
    \end{gather*}
\end{assumption}
\begin{remark}\label{remark:F_expression}
In Assumption \ref{assmp:affine_F}, suppose that
$$
    M(\xi)=M\cdot \mathrm{repvec}(\xi;n)+M_0,\quad q(\xi)=Q\xi+q_0,
$$
where
\begin{gather*}
    M:=\left[\begin{array}{cccc}
        (m^{1,1})^\top & (m^{1,2})^\top & \dots & (m^{1,n})^\top \\
        (m^{2,1})^\top & (m^{2,2})^\top & \dots & (m^{2,n})^\top \\
        \vdots & \vdots & \ddots & \vdots \\
        (m^{n,1})^\top & \dots & \dots & (m^{n,n})^\top
    \end{array}\right]\in\Re^{n\times mn},\\
    \mathrm{repvec}(\xi;n):=\left[\begin{array}{cccc}
        \xi & & & \\
        & \xi & & \\
        & & \ddots & \\
        & & & \xi
    \end{array}\right]\in\Re^{mn\times n},\quad
    M_0:=[\bar{m}^1_0,\bar{m}^2_0,\dots,\bar{m}^n_0]^\top\in\Re^{n\times n},\\
    Q:=[q^1,q^2,\dots,q^n]^\top\in\Re^{n\times m},\ q_0:=[q^1_0,q^2_0,\dots,q^n_0]\in\Re^n.
\end{gather*}
Then, $F(x,\xi)$ can also be written as
$$
    F(x,\xi)=\left(M\cdot\mathrm{repvec}(\xi;n)+M_0\right)x+(Q\xi+q_0).
$$
\end{remark}

Let us introduce the convexity of nonlinear matrix-valued functions and its related property.
\begin{definition}[{Shapiro \cite{Shapiro1997}}]\label{def:matrix_convexity}
	A nonlinear matrix-valued function $X\colon\Re^m\rightarrow\mathbb{S}^n$ is said to be positive semidefinite (psd-) convex if
	\begin{align}\label{ieq:psd-convex}
		X(\gamma x + (1-\gamma) y) - \gamma X(x) - (1-\gamma) X(y) \preceq O
	\end{align}
	for all $x,y\in\Re^m$ and $\gamma\in[0,1]$.
\end{definition}

\begin{proposition}\label{prop:psd-convexity}
	The mapping $X$ is psd-convex if and only if for any $v\in\Re^n$ with $v_1=1$, the function $\phi(\cdot;v)\colon\Re^m\to\Re$ defined by
	\[
	    \phi(x;v)\coloneqq \Square{1,\ v^\top}X(x)\left[\begin{array}{c} 1 \\ v\end{array}\right]
	\]
	is convex with respect to $x\in\Re^m$.
\end{proposition}
\proof
Lemma~\ref{lem:NSD_constraint} ensures that matrix inequality \eqref{ieq:psd-convex} is equivalent to
\[
    \Square{1,\ v'^\top}
    (X(\gamma x+(1-\gamma)y)-\gamma X(x)-(1-\gamma)X(y))
    \Square{\begin{array}{c}1\\v'\end{array}}\leq 0,
\]
for any $v'\in\Re^{n-1}$.
Hence, we have
\[
    \phi(\gamma x+(1-\gamma)y;v)\leq\gamma\phi(x;v)+(1-\gamma)\phi(y;v)
\]
for any $v\in\Re^n$ with $v_1=1$.
Therefore, $X$ is psd-convex if and only if $\phi(\cdot,v)$ is convex with respect to $x\in\Re^m$ for every $v\in\Re^n$.
\qed

\bigskip

We show the convexity of (NSDP) and $\mathrm{(NSDP')}$.

\begin{theorem}\label{thm:nsd_concavity_of_D}
    Suppose that Assumption~{\rm\ref{assmp:affine_F}} holds and that the matrix $M(\xi)$ defined in Assumption~{\rm\ref{assmp:affine_F}} satisfies the following condition: There exists $\beta_0>0$ such that
    \begin{align}\label{ieq:uniformly_pd}
        \inf_{\xi\in\Xi,\|v\|=1} v^\top M(\xi) v \geq \beta_0.
    \end{align}
    Then the matrix-valued function $-\D$ is psd-convex for all $\alpha\geq 1/(2\beta_0)${\rm;} thus, {\rm(NSDP)} and $\mathrm{(NSDP')}$ are convex.
\end{theorem}
\proof
Note that if the matrix-valued function $-\D$ is psd-convex, (NSDP) and $\mathrm{(NSDP')}$ are convex optimization problems.
Therefore, we verify that $-\D$ is psd-convex for all $\alpha\geq 1/(2\beta_0)$.

Suppose that $\alpha\geq 1/(2\beta_0)$.
Proposition~\ref{prop:psd-convexity} states that $-\D$ is psd-convex if and only if for all $\xi\in\Re^m$, the following function $\phi_\alpha(\cdot,\xi)\colon\mathcal{W}\rightarrow\Re$ is convex with respect to $w$:
\begin{align*}
    &\phi_\alpha(w;\xi)\coloneqq\Square{1, \xi^\top}(-\D(w))\Square{\begin{array}{c} 1 \\ \xi\end{array}}= -y_0-\xi^\top y - \xi^\top Y\xi + \omega_\alpha(x,\lambda,\mu;\xi),
\end{align*}
where the last equality follows from \eqref{D_alpha=omega}.

Now, since the function $\phi_\alpha(\cdot,\xi)$ is linear with respect to $(y_0,y,Y)$, it suffices to show that $\omega_\alpha$ is convex with respect to $(x,\lambda,\mu)$ for all $\alpha\geq1/(2\beta_0)$.
The Hessian of $\omega_\alpha$ in regard to $(x,\lambda,\mu)$ is given by
\[
    \small
    \nabla^2_{(x,\lambda,\mu)}\omega_\alpha(x,\lambda,\mu;\xi)=\alpha \Square{
        \begin{array}{ccc}
            M(\xi)^\top M(\xi) & (M(\xi)-\frac{1}{\alpha}I)^\top A^\top & -M(\xi)^\top + \frac{1}{\alpha}I \\
            A(M(\xi)-\frac{1}{\alpha}I) & AA^\top & -A \\
            -M(\xi)+\frac{1}{\alpha}I & -A^\top & I
        \end{array}
    }
    \normalsize.
\]
By considering the Schur complement of the above matrix,
\small
\begin{align*}
    &\Square{
        \begin{array}{cc}
            M(\xi)^\top M(\xi) & (M(\xi)-\frac{1}{\alpha}I)^\top A^\top \\
            A(M(\xi)-\frac{1}{\alpha}I) & AA^\top
        \end{array}
    }- \Square{
        \begin{array}{c}
            -M(\xi)^\top + \frac{1}{\alpha} I \\ -A
        \end{array}
    }\Square{-M(\xi)+\frac{1}{\alpha} I\quad -A^\top} \\
    & = \frac{1}{\alpha}\Square{
        \begin{array}{cc}
            (M(\xi)^\top + M(\xi))-\frac{1}{\alpha} I & O \\
            O & O
        \end{array}
    }\succeq O
\end{align*}
\normalsize
if and only if $\nabla^2_{(x,\lambda,\mu)}\omega_\alpha(x,\lambda,\mu;\xi)\succeq O$.
Since $(M(\xi)^\top + M(\xi))-1/\alpha\ I \succeq O$ from $\alpha\geq1/(2\beta_0)$, it can be easily seen that $\nabla^2_{(x,\lambda,\mu)}\omega_\alpha(x,\lambda,\mu;\xi)\succeq O$, i.e., $-\D$ is psd-convex for all $\alpha\geq1/(2\beta_0)$. Hence, (NSDP) and $\mathrm{(NSDP')}$ are convex optimization problems.
\qed


\begin{remark}\label{remark:convex_regularizer}
    Condition \eqref{ieq:uniformly_pd} is rather restrictive for some applications.
    One remedy is to add a proximal term $\epsilon(x-x^k)$ to the mapping $F$, where $\epsilon>0$ is a sufficiently small constant.
\end{remark}

\begin{remark}\label{remark:Zhu_convexity}
    When $S=\Re^n_+$, problem \eqref{prob:ZhuDRCP} for the SLCP proposed by Zhu et al. \cite{Zhu2017} may not be reformulated as a convex NSDP because the objective function $\Psi(x,\xi)=\|x\circ F(x,\xi)\|^2_2$ is not convex with respect to $x$ in general.
\end{remark}

Although we adopt the regularized gap function for the NSDP approximation, the similar results may also be obtained by utilizing another merit function, such as $f_\infty(x,\xi)\coloneqq\max_{z\in S}\langle F(x,\xi),x-z\rangle$.
However, it would be necessary to discuss whether the DRERM with $f_\infty$ is reasonable method for solving the SVIP.
In fact, the ERM with $f_\infty$ may be unsuitable to measure the distance to solutions of SVIP \eqref{SVIP} because $f_\infty(x,\xi)$ takes $+\infty$ for some $x\in S$ and is not differentiable in general.
For such reasons, we did not adopt $f_\infty$ for (DRERM).

\section{Numerical experiments}\label{sec:numerical_experiments}

This section provides numerical results to demonstrate the validity of the DRERM model.
In particular, we first compare the DRERM with the ERM proposed by Luo and Lin \cite{Luo2009a} in terms of robustness.
Second, we quantitatively investigate the robustness of solutions obtained from the DRERM model when the confidence parameters $\gamma_1$ and $\gamma_2$ for the mean and variance of the ambiguity set $\mathscr{P}$, respectively, are gradually changed.

Throughout this section, we use the following example.
\begin{example}[Two-person noncooperative games]
	Two players are competing with each other to minimize their own cost functions.
	Each player $\nu\in\{1,2\}$ solves the following optimization problem:
	\begin{align}\label{ex:prob:player}
		\begin{array}{cl}
			\underset{x^\nu\in\Re^{n_\nu}}{\min} & \displaystyle\frac{1}{2}\Paren{x^\nu}^\top M_\nu x^\nu + v^\nu(x^\nu,x^{-\nu},\xi) + q^\nu(\xi)^\top x^\nu \\
			\text{s.t.} & A_\nu x^\nu\leq b^\nu,
		\end{array}
	\end{align}
	where $M_\nu\in\mathbb{S}^{n_\nu}_{++}$, $A_\nu\in\Re^{l_\nu\times n_\nu}$, $b^\nu\in\Re^{l_\nu}$, and $q^\nu(\xi)\in\Re^{n_\nu}$.
	Here, $v^\nu(x^\nu,x^{-\nu},\xi)$ is a zero-sum function defined by
	\begin{align*}
		v^\nu(x^\nu,x^{-\nu},\xi)\coloneqq
		\begin{cases*}
			\begin{array}{lc}
				\Paren{x^1}^\top R(\xi) x^2 & \text{if } \nu=1, \\
				-\Paren{x^2}^\top R(\xi)^\top x^1 & \text{if } \nu=2,
			\end{array}
		\end{cases*}
	\end{align*}
	where $R(\xi)\in\Re^{n_1\times n_2}$, and $x^{-\nu}\in\Re^{n_{-\nu}}$ is the decision variable of the rival player. \qed
\end{example}

The above noncooperative game can be reformulated as SVIP \eqref{SVIP} when the mapping $F(\cdot,\xi)\colon\Re^n\to\Re^n$ and the set $S\subset\Re^n$ are given as follows:
\begin{align}
	F(x,\xi)&=\Square{
		\begin{array}{cc}
			M_1 & R(\xi) \\
			-R(\xi)^\top & M_2
		\end{array}
		}x + \Square{
		\begin{array}{c}
			q^1(\xi) \\
			q^2(\xi)
		\end{array}
	},\label{eq:mapping_F}\\
	S&= \Curly{x\in\Re^n\ \middle|\ \Square{
			\begin{array}{cc}
				A_1 & O \\
				O & A_2
			\end{array}
			}x \leq \Square{
			\begin{array}{c}
				b^1 \\
				b^2
			\end{array}
	}},\nonumber\\
	x&=\Square{(x^1)^\top,\ (x^2)^\top}^\top\in\Re^{n_1+n_2}. \nonumber
\end{align}
Note that it is easy to verify that the coefficient matrix in \eqref{eq:mapping_F} satisfies the assumption of Theorem~\ref{thm:nsd_concavity_of_D}; hence, we solve a convex NSDP in the experiments.

We generate numerical instances of problem \eqref{ex:prob:player} according to the following manners:
\begin{itemize}
\item We set $n_1=n_2=2$, $m=n_1n_2+2=6$, and $l_1=l_2=2$.
\item The matrix $M_\nu$ is generated by $L_\nu {L_\nu}^\top + I$, where the matrix $L_\nu\in\Re^{2\times 2}$ is lower triangular and its elements are randomly generated from the interval $[-5,5)$.
\item Each element of the matrix $A_\nu\in\Re^{2\times 2}$ and the vector $b^\nu\in\Re^2$ is randomly generated from $[-2,2)$ and $[0, 10)$, respectively.
\item We set the regularization parameter $\alpha$ by $1/\beta_0$ to ensure that the derived NSDP is convex, where $\beta_0$ is the minimum eigenvalue of the matrix
	\[
		\Square{\begin{array}{cc}
				M_1 & O_{n_1\times n_2} \\
				O_{n_2\times n_1} & M_2
			\end{array}            
		}\in\Re^{4 \times 4}.
	\]
\item We define the random variable $\xi\in\Re^m$ by $\xi=[\xi_1,\dots,\xi_6]^\top$.
\item The matrix $R(\xi)$ is defined by
    $$
        R(\xi)\coloneqq\left[
        \begin{array}{cc}
            \xi_1 & \xi_2 \\
            \xi_3 & \xi_4
        \end{array}
        \right]+R_0\in\Re^{2\times 2},\ 
        R_0\coloneqq\Square{
			\begin{array}{ccc}
				r^{1,1}_0 & r^{1,2}_0 \\
				r^{2,1}_0 & r^{2,2}_0
			\end{array}
		}\in\Re^{2\times 2}
    $$
    where $r_0^{i,j}$, $i,j=1,2$ are nominal values generated randomly from $[-5,5)$.
\item The vector $q(\xi)\coloneqq\Paren{q^1(\xi)^\top,q^2(\xi)^\top}^\top\in\Re^4$ is defined by
$$q(\xi)=Q\xi+q_0,$$ where
	\begin{gather*}
		{{Q}}=\left[
			\begin{array}{cccccc}
			    0 & 0 & 0 & 0 & 1 & 0 \\
			    0 & 0 & 0 & 0 & 1 & 0 \\
			    0 & 0 & 0 & 0 & 0 & 1 \\
			    0 & 0 & 0 & 0 & 0 & 1
			\end{array}
		\right],\ 
		q_0=-\Square{
		    \begin{array}{cc}
   				M_1 & R_0 \\
   				-R_0^\top & M_2
   			\end{array}            
	    }x^*_0, 
   	\end{gather*}
	and the vector $x^*_0\in\Re^{4}$ is randomly generated from $[-2,2)$.
\end{itemize}

In the experiments, all programs are implemented with Python 3.8 and run on a machine with Intel Core i7-8700K @ 3.70GHz CPU and 32 GB RAM.

\subsection{Comparison to the ERM model}\label{sec:numerical_experiments:comp_ERM}

Here, we suppose that $\Xi=\Re^6$ and $\xi$ follows the normal distribution $\mathcal{N}(\mu_0,\Sigma_0)$, where the mean $\mu_0$ and the variance-covariance matrix $\Sigma_0$ are given as follows:
\begin{align}\label{eq:Sigma0}
    \mu_0 = 0,\ 
    \Sigma_0=\Square{
       \begin{array}{cccc}
           2 & 1.6 & \cdots & 1.6 \\
           1.6 & 2 & \cdots & 1.6 \\
           \vdots & \vdots & \ddots & \vdots \\
           1.6 & 1.6 & \cdots & 2
       \end{array}
   }.
\end{align}

In the ERM model, we use the regularized gap function $f_\alpha$ proposed by Luo and Lin \cite{Luo2009a} as the merit function $f$.
In the experiments, because it is difficult to exactly compute the expected value $\mathbb{E}[f_\alpha(x,\xi)]$, we obtain its approximate value using a quasi-Monte Carlo method described below:
\[
    \mathbb{E}[f_\alpha(x,\xi)]\approx \theta^k(x)\coloneqq\displaystyle\frac{1}{N_k}\Sum{\hat{\xi}^k\in\Xi^k}{}{f_\alpha(x,\hat{\xi}^k)p(\hat{\xi}^k)},
\]
where the uniform random vector $\hat{\xi}^k\in\Xi^k$ is generated by 
\[
    \hat{\xi}^k=\Paren{(\mu_0-3\sqrt{2})+(\mu_0+3\sqrt{2})\zeta^i}\mathbf{1}_m,
\]
and $\zeta^i$ is a Sobol point from the interval $[0,1)$.
The set $\Xi^k\coloneqq\{\hat{\xi}^i\mid i=1,2,\ldots,N_k\}\subset\Xi$ is the collection of the samples $\hat{\xi}^k$, which approximates the support $\Xi$, and $p(\cdot)$ is the probability density function of the normal distribution $\mathcal{N}(\mu_0,\Sigma_0)$.
Note that as the number of samples $N_k$ and dimensions $m$ increased, it may face underflow and subsequently fail to evaluate $\theta^k(x)$.
To avoid this, we multiply $\theta^k(x)$ by $1/p(\mu_0)$.
Summarizing the above arguments, we solve the following approximate problem for (ERM) with the regularized gap function:
\begin{align}\label{prob:approx_ERM}
    \begin{array}{cl}
        \text{min} & \theta^k(x)/p(\mu_0)\\
        \text{s.t.} & x\in S,
    \end{array}
\end{align}
We use SLSQP package, which is based on sequential quadratic programming methods, in Scipy.Optimize module to obtain a solution to problem \eqref{prob:approx_ERM}.
The initial point is set to $0$, and the termination criterion for the residual of the Karush--Kuhn--Tucker condition is set to $10^{-7}$.

In the DRERM, because we know the exact values $\mu_0$ and $\Sigma_0$ in advance, the ambiguity set $\mathscr{P}$ is given by \eqref{ambiguity0}.
When $\Xi=\Re^6$ and $\mathscr{P}$ is given as \eqref{ambiguity0}, (DRERM) can be reformulated as the following NSDP, which can be regarded as the special case of $\mathrm{(NSDP')}$:
\begin{align}\label{prob:NSDP_specialized}
    \begin{array}{cl}
        \underset{(x,\lambda,y_0,y,Y)}{\text{min}} & y_0 + \mu_0^\top y + \Ip{\Sigma_0+\mu_0\mu_0^\top}{Y} \\
        \text{s.t.} & \mathcal{D}_\alpha(x,\lambda,y_0,y,Y) \succeq O, \\
        & Ax\leq b,\ \lambda\in\Re^2_-,
    \end{array}
\end{align}
where $\Re^2_-\coloneqq\{\lambda\in\Re^2\mid \lambda\leq 0\}$.
To solve \eqref{prob:NSDP_specialized}, we utilize an interior point method, which is a hybrid method of \cite{Yamashita2009} and \cite{Yamashita2012}.
The initial point and termination criterion are the same as the method for \eqref{prob:approx_ERM}.

We prepare 10 numerical instances of SVIP~\eqref{SVIP} and solve them via \eqref{prob:approx_ERM} and \eqref{prob:NSDP_specialized}, where we set two cases where $N_k=80$ and $N_k=10000$ in \eqref{prob:approx_ERM}.
Let $x^{i*}_{\mathrm{ERM}}$ and $x^{i*}_{\mathrm{DRERM}}$ be solutions to \eqref{prob:approx_ERM} and \eqref{prob:NSDP_specialized} at the $i$-th instance, respectively.
In what follows, for a realization $\bar{\xi}^j$ of the random variable $\xi$, $\realgap{ERM}$ and $\realgap{DRERM}$ respectively denote $f_\alpha(x^{i*}_{\mathrm{ERM}},\bar{\xi}^j)$ and $f_\alpha(x^{i*}_{\mathrm{DRERM}},\bar{\xi}^j)$ for simplicity.

To quantitatively evaluate the solutions $x^{i*}_{\mathrm{ERM}}$ and $x^{i*}_{\mathrm{DRERM}}$, we conduct the following steps:
\begin{enumerate}
    \item[(i)] Generate $N\coloneqq$ 5000 realizations $\{\bar{\xi}^j\}_{j=1}^N$, where each realization $\bar{\xi}^j$ follows the normal distribution $\mathcal{N}(\mu_1,\Sigma_1)$.
    Here, $\mu_1$ and $\Sigma_1$ are respectively the perturbations of $\mu_0$ and $\Sigma_0$ as follows:
    $$
        \mu_1\coloneqq\mu_0+\delta_\mu,\ \Sigma_1\coloneqq\Sigma_0+\Delta_\Sigma,
    $$
    where each element of $\delta_\mu\in\Re^6$ and $\Delta_\Sigma\in\mathbb{S}^6$ are uniformly generated from the interval $[-0.1, 0.1)$.
    \item[(ii)] Compute the regularized gap function values $\{\realgap{ERM}\}_{j=1}^{N}$ and $\{\realgap{DRERM}\}_{j=1}^{N}$ by using the realizations $\{\bar{\xi}^j\}_{j=1}^N$ for each solution.
    \item[(iii)] Evaluate the solutions $x^{i*}_{\mathrm{ERM}}$ and $x^{i*}_{\mathrm{DRERM}}$ by using the following five indicators, which represent the rates of change (RC):
\begin{itemize}
	\item Minimum:
		\begin{align}\label{eq:difference_of_min}
			(\min_j \realgap{DRERM} - \min_j \realgap{ERM})/\min_j \realgap{ERM},
		\end{align}
	\item Maximum:
		\begin{align}\label{eq:difference_of_max}
			(\max_j \realgap{DRERM} - \max_j \realgap{ERM})/\max_j \realgap{ERM},
		\end{align}
	\item Mean:
		\begin{align}\label{eq:difference_of_mean}
			(\meangap{DRERM} - \meangap{ERM})/\meangap{ERM},
		\end{align}
		where $\meangap{\cdot}\coloneqq\frac{1}{N}\sum_{j=1}^{N}\realgap{\cdot}$.
	\item Median:
		\begin{align}\label{eq:difference_of_median}
			(\medgap{DRERM} - \medgap{ERM})/\medgap{ERM},
		\end{align}
		where $\medgap{\cdot}\coloneqq(f^{i[N/2]}_{\cdot}+f^{i[N/2+1]}_\cdot)/2$, and $f^{i[j]}_{\cdot}$ denotes the $j$-th largest regularized gap function value in the 5000 realizations.
	\item Standard deviation (SD):
	\begin{align}\label{eq:difference_of_SD}
		(\vargap{DRERM}-\vargap{ERM})/\vargap{ERM},
	\end{align}
	where $\vargap{\cdot}\coloneqq\sqrt{\frac{1}{N-1}\sum_{j=1}^{N}(\realgap{\cdot}-\meangap{\cdot})^2}$.
\end{itemize}
\end{enumerate}

\begin{figure}[htb]
    \centering
    \subfloat[RC \eqref{eq:difference_of_min} when $N_k=$ 80]{%
        \resizebox*{7cm}{!}{\includegraphics{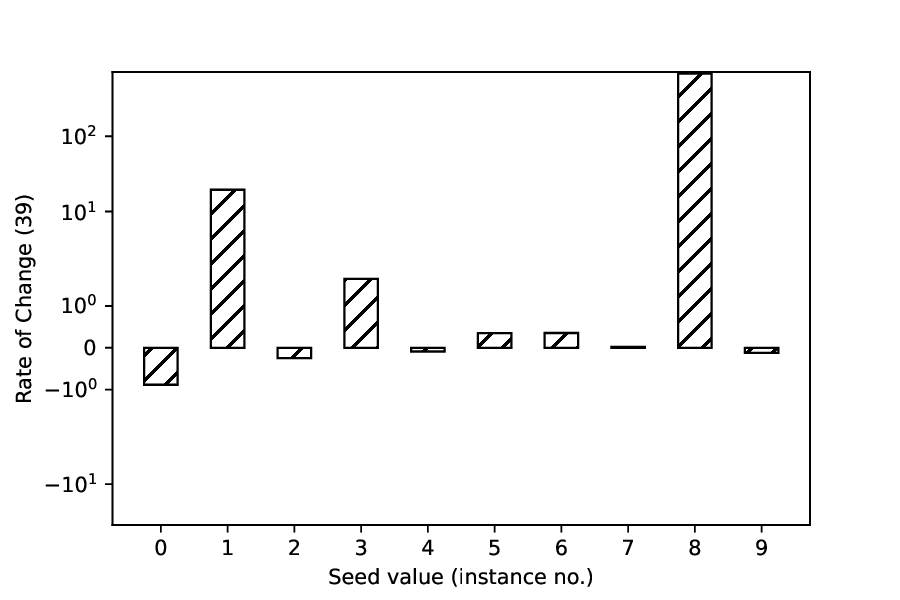}}\label{fig:ERM_comp:min_QMC80}
    }
    \subfloat[RC \eqref{eq:difference_of_min} when $N_k=$ 10000]{%
        \resizebox*{7cm}{!}{\includegraphics{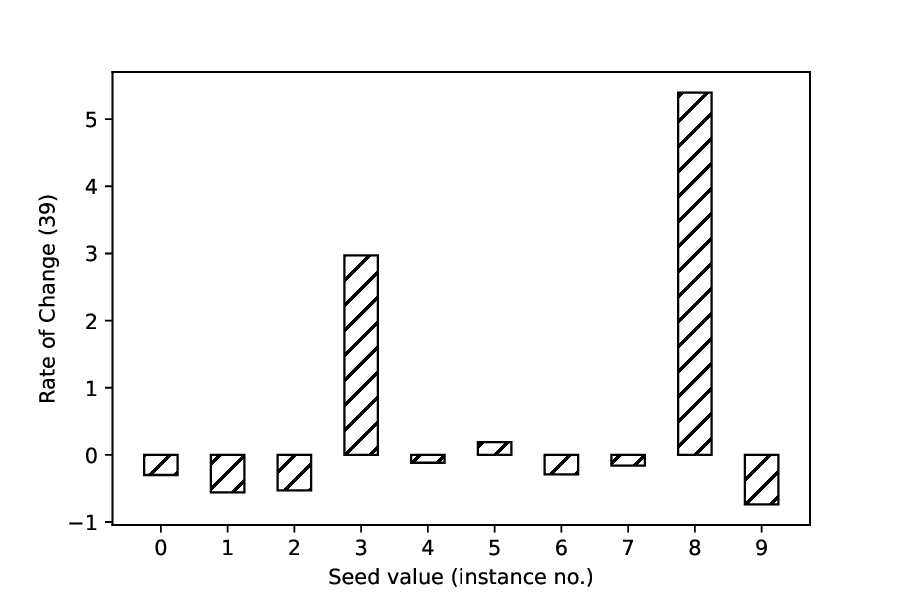}}\label{fig:ERM_comp:min_QMC10000}
    }\\
    \subfloat[RC \eqref{eq:difference_of_max}--\eqref{eq:difference_of_SD} when $N_k=$ 80]{%
        \resizebox*{7cm}{!}{\includegraphics{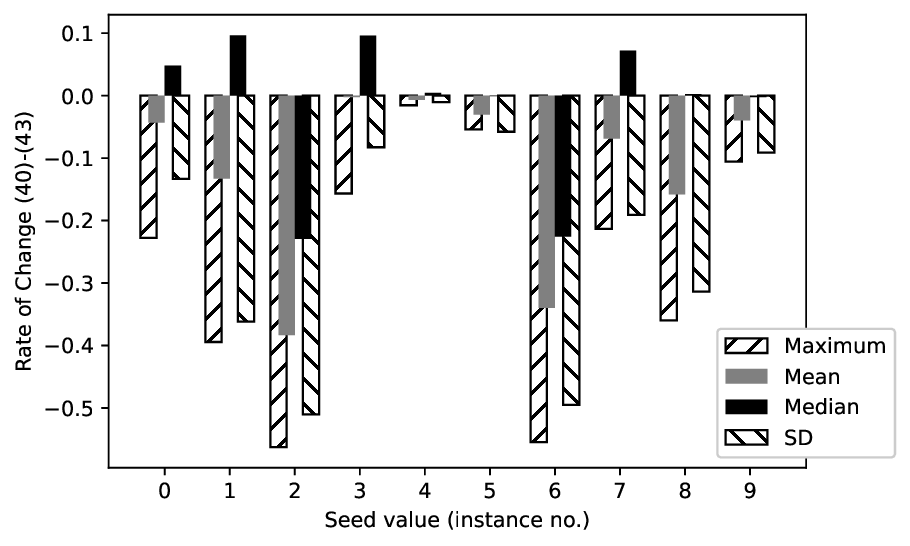}}\label{fig:ERM_comp:others_QMC80}
    }
    \subfloat[RC \eqref{eq:difference_of_max}--\eqref{eq:difference_of_SD} when $N_k=$ 10000]{%
        \resizebox*{7cm}{!}{\includegraphics{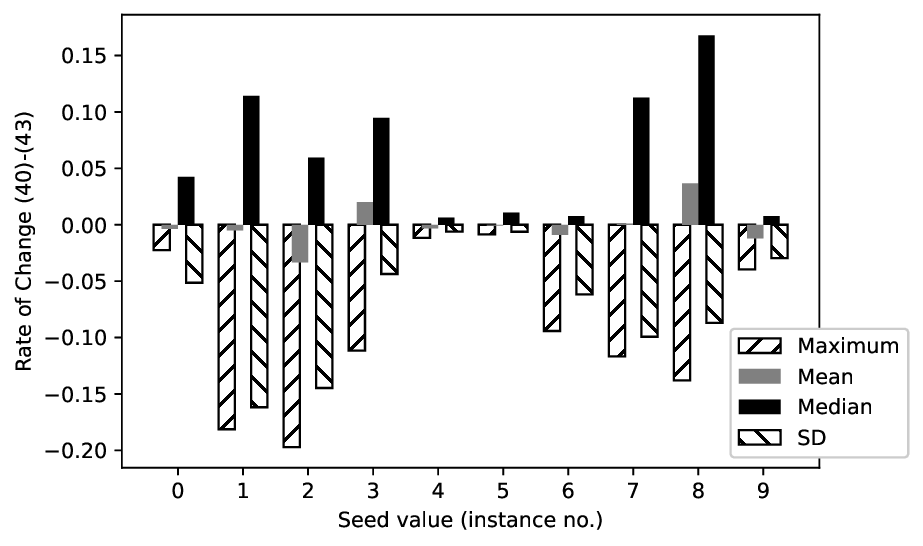}}\label{fig:ERM_comp:others_QMC10000}
    }
    \caption{The rate of change between the ERM and the DRERM.}\label{fig:ERM_comp}
\end{figure}

The computational results are shown in Figure~\ref{fig:ERM_comp}.
In each graph, the horizontal and the vertical axes represent the instance number and the RC, respectively.
Figures~\ref{fig:ERM_comp}\subref{fig:ERM_comp:min_QMC80} and~\ref{fig:ERM_comp}\subref{fig:ERM_comp:min_QMC10000} indicate the RC evaluated by \eqref{eq:difference_of_min} for $N_k=80$ and $N_k=10000$, respectively, and Figures~\ref{fig:ERM_comp}\subref{fig:ERM_comp:others_QMC80} and~\ref{fig:ERM_comp}\subref{fig:ERM_comp:others_QMC10000} represent the RC evaluated by \eqref{eq:difference_of_max}--\eqref{eq:difference_of_SD} for each $N_k$.
Note that the vertical axis of Figure~\ref{fig:ERM_comp}\subref{fig:ERM_comp:min_QMC80} is a logarithmic scale.

First, we focus on the minimum values, i.e., Figures~\ref{fig:ERM_comp}\subref{fig:ERM_comp:min_QMC80} and~\ref{fig:ERM_comp}\subref{fig:ERM_comp:min_QMC10000}.
We observe that for most of instances of $N_k=80$ and $N_k=10000$, the minimum values of the ERM tend to be small compared with the DRERM.
In particular, the 8-th instance in Figure~\ref{fig:ERM_comp}\subref{fig:ERM_comp:min_QMC80} indicates a significant difference between the ERM and DRERM models.
Indeed, $\min_j f^{8j}_{\mathrm{ERM}}=0.0023$ and $\min_j f^{8j}_{\mathrm{DRERM}}=1.5784$, and they have a 690-fold difference.
In the case of $N_k=$ 10000, the gaps between the ERM and DRERM are small for all instances compared with $N_k=80$.

Next, we focus on Figures~\ref{fig:ERM_comp}\subref{fig:ERM_comp:others_QMC80} and~\ref{fig:ERM_comp}\subref{fig:ERM_comp:others_QMC10000}.
Notably, the values of the gap function of maximum and SD on the DRERM are smaller than the ERM for all instances for $N_k=80$ and $N_k=10000$.
This is an important result that shows that the DRERM is reasonably designed to consider the distributionally worst case in terms of the expected value of the regularized gap function.

From the above results, we confirm that the DRERM can obtain more robust solutions that consider outliers, while the ERM is not as robust as the DRERM even when $N_k$ is sufficiently large in spite of using the exact distribution function for evaluating the expected value.
This is because the ERM is designed to minimize the expected value of the regularized gap function; hence, it cannot directly consider the variance and maximum value.
In fact, the median of $\realgap{ERM}$ with $N_k=$ 10000 is less than the DRERM; however, outliers of realizations $\bar{\xi}^j$ adversely affect the mean of the regularized gap values.
As a result, the difference between the mean of $\realgap{ERM}$ with $N_k=$ 10000 and that of $\realgap{DRERM}$ is insignificant.


\subsection{Analysis of solution by varying confidence parameters}\label{sec:numerical_experiments:drerm_bench}

In this section, we assume that $\Xi=\Re^6$ and the estimated mean $\tilde{\mu}_0$ and variance-covariance matrix $\tilde{\Sigma}_0$ are given as follows:
\[
	\tilde{\mu}_0 \coloneqq \mu_0 + u^6,\ \tilde{\Sigma}_0 \coloneqq \Sigma_0 + U_6,
\]
where each element of $u^6\in\Re^6$ and $U_6\in\mathbb{S}^6$ are uniformly generated from $[-0.25, 0.25)$ and $[-0.2, 0.2)$, respectively.
Here, the true $\mu_0$ and $\Sigma_0$ are the same as \eqref{eq:Sigma0}, and the confidence regions of $\tilde{\mu}_0$ and $\tilde{\Sigma}_0$ in the ambiguity set $\mathscr{P}$ are given as follows:
\begin{align}
	\Paren{\ExpP{\xi}-\tilde{\mu}_0}^\top\tilde{\Sigma}_0^{-1}\Paren{\ExpP{\xi}-\tilde{\mu}_0}\leq \gamma_1,\label{ex:ambiguity_set1}\\
	\ExpP{\Paren{\xi-\tilde{\mu}_0}\Paren{\xi-\tilde{\mu}_0}^\top}\preceq\gamma_2\tilde{\Sigma}_0.\label{ex:ambiguity_set2}
\end{align}
In this setting, we solve the following NSDP:
\begin{align}\label{prob:equivalent_NSDP2}
	\begin{array}{cl}
		\underset{(x,\lambda,y_0,y,Y,z_0)}{\min} & z_0 + y_0 + \tilde{\mu}_0^\top y + \langle \gamma_2\tilde{\Sigma}_0+\tilde{\mu}_0\tilde{\mu}_0^\top, Y\rangle \\
		\text{s.t.} & z_0 \geq \sqrt{\gamma_1}\left\|\tilde{\Sigma}_0^{1/2}(y+2Y\tilde{\mu}_0)\right\|, \\
								& \displaystyle\mathcal{D}_\alpha(x,\lambda,y_0,y,Y)\succeq O, \\
								& Ax\leq b,\ \lambda\in\Re^2_-.
	\end{array}
\end{align}
Here, we solve~\eqref{prob:equivalent_NSDP2} using the interior point method, which is the same method for solving~\eqref{prob:NSDP_specialized}.
The initial point is set as 0, and the stopping criterion is $10^{-7}$.
Note that we set $\alpha>0$ to ensure that problem \eqref{prob:equivalent_NSDP2} is convex.
Let $x^*_{\gamma_1,\gamma_2}$ be a solution of problem \eqref{prob:equivalent_NSDP2} for given $\gamma_1$ and $\gamma_2$.

In the first experiment, we quantitatively analyze the characteristics of the solutions in the case where $\gamma_1$ is incremented by 0.1 from 0.1 to 2, and $\gamma_2$ is set to 1 or 2.
We prepare realizations $\{\bar{\xi}^j\}_{j=1}^N$, where each $\bar{\xi}^j$ follows $\mathcal{N}(\mu_0,\Sigma_0)$ and $N=$ 5000.
After obtaining a solution $x^*_{\gamma_1,\gamma_2}$, we compute the maximum, mean, and SD of $\{f_\alpha(x^*_{\gamma_1,\gamma_2},\bar{\xi}^j)\}_{j=1}^N$.

\begin{figure}[htbt]
    \centering
    \subfloat[Maximum for fixed $\gamma_2=1$]{%
        \resizebox*{7cm}{!}{\includegraphics{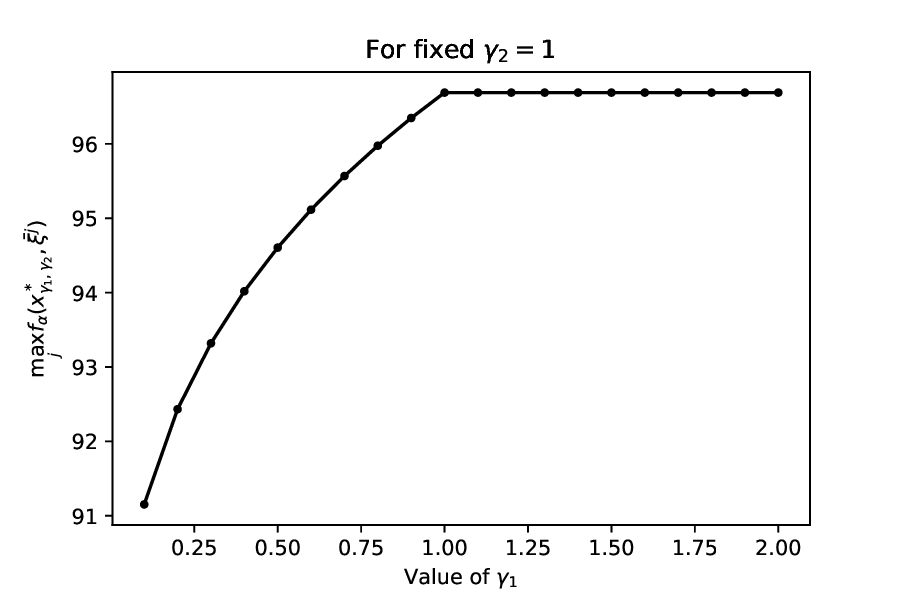}}\label{fig:DRERM_Bench_gamma1:gamma_2=1_max}
    }
    \subfloat[Maximum for fixed $\gamma_2=2$]{%
        \resizebox*{7cm}{!}{\includegraphics{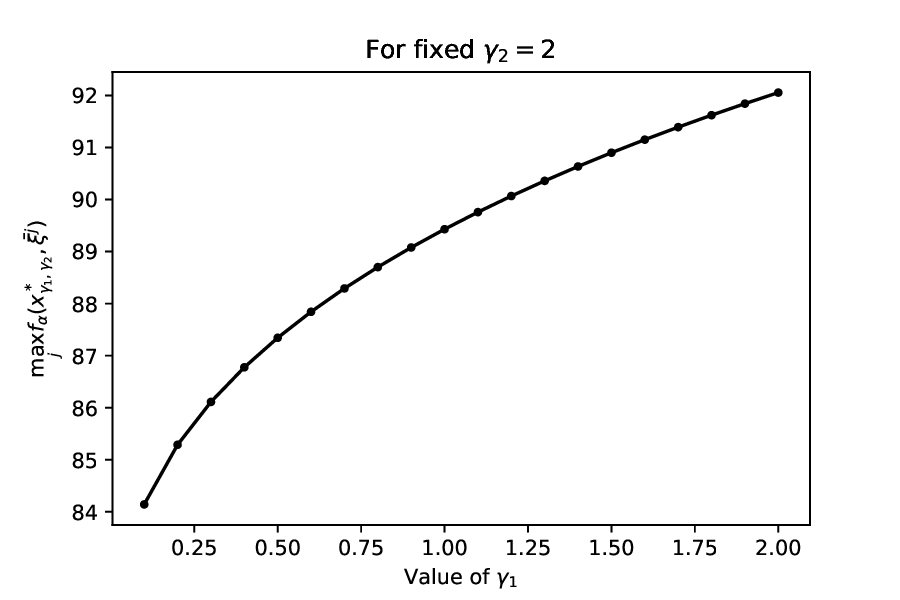}}\label{fig:DRERM_Bench_gamma1:gamma_2=2_max}
    }\\
    \subfloat[Mean and SD for fixed $\gamma_2=1$]{%
        \resizebox*{7cm}{!}{\includegraphics{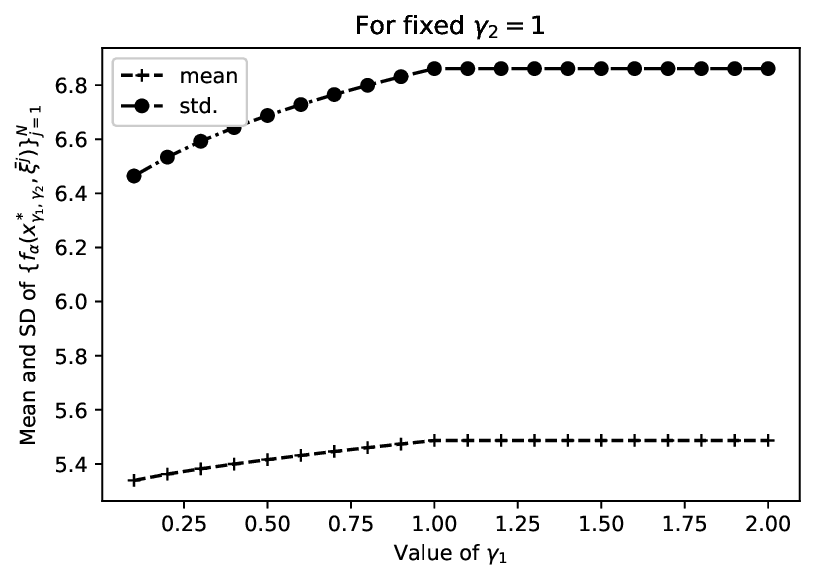}}\label{fig:DRERM_Bench_gamma1:gamma_2=1_others}
    }
    \subfloat[Mean and SD for fixed $\gamma_2=2$]{%
        \resizebox*{7cm}{!}{\includegraphics{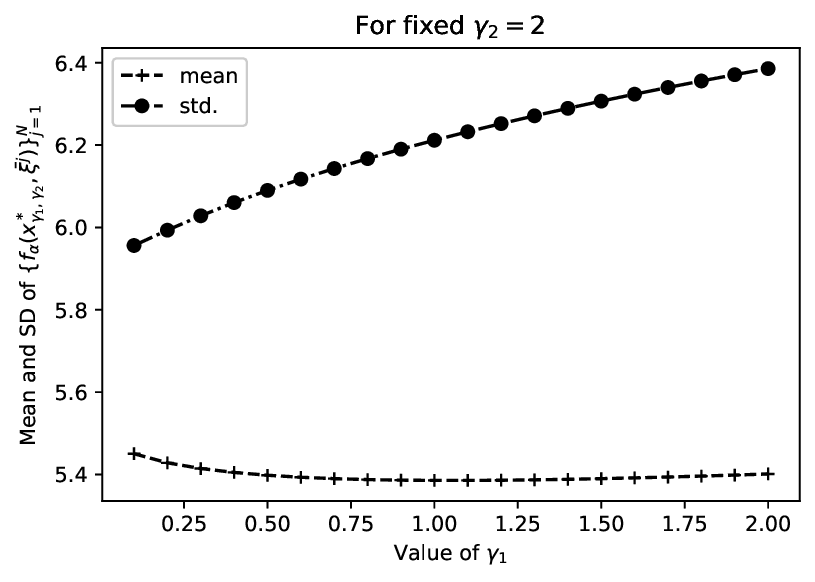}}\label{fig:DRERM_Bench_gamma1:gamma_2=2_others}
    }
    \caption{Maximum, mean, and SD of 5000 realizations of the regularized gap function when $\gamma_1$ is varied.}\label{fig:DRERM_Bench_gamma1}
\end{figure}

Figure~\ref{fig:DRERM_Bench_gamma1} shows the results of the first experiment.
In each graph, the horizontal and vertical axes represent the values of $\gamma_1$ and the regularized gap function, respectively.
The curves in Figures~\ref{fig:DRERM_Bench_gamma1}\subref{fig:DRERM_Bench_gamma1:gamma_2=1_max} and~\ref{fig:DRERM_Bench_gamma1}\subref{fig:DRERM_Bench_gamma1:gamma_2=2_max} indicate the maximum of $\{f_\alpha(x^*_{\gamma_1,\gamma_2},\bar{\xi}^j)\}_{j=1}^N$ for fixed $\gamma_2=1$ and $\gamma_2=2$, respectively, and Figures~\ref{fig:DRERM_Bench_gamma1}\subref{fig:DRERM_Bench_gamma1:gamma_2=1_others} and~\ref{fig:DRERM_Bench_gamma1}\subref{fig:DRERM_Bench_gamma1:gamma_2=2_others} represent the mean and SD of $\{f_\alpha(x^*_{\gamma_1,\gamma_2},\bar{\xi}^j)\}_{j=1}^N$ for fixed $\gamma_2=1$ and $\gamma_2=2$, respectively.

In Figures~\ref{fig:DRERM_Bench_gamma1}\subref{fig:DRERM_Bench_gamma1:gamma_2=1_max} and~\ref{fig:DRERM_Bench_gamma1}\subref{fig:DRERM_Bench_gamma1:gamma_2=1_others} (when $\gamma_2=1$), the maximum, mean, and SD of regularized gap function values increase as $\gamma_1$ increases.
However, Figures~\ref{fig:DRERM_Bench_gamma1}\subref{fig:DRERM_Bench_gamma1:gamma_2=2_max} and~\ref{fig:DRERM_Bench_gamma1}\subref{fig:DRERM_Bench_gamma1:gamma_2=2_others} (when $\gamma_2=2$) indicate that the values of the maximum and SD are entirely smaller than the case where $\gamma_2=1$; we will discuss the reason in the next experiment.
In particular, from Figure~\ref{fig:DRERM_Bench_gamma1}\subref{fig:DRERM_Bench_gamma1:gamma_2=2_others}, the curve of the mean gradually decreases for $0.1\leq\gamma_1\leq 1$, unlike the case where $\gamma_2=1$.
Moreover, Figures~\ref{fig:DRERM_Bench_gamma1}\subref{fig:DRERM_Bench_gamma1:gamma_2=1_max} and~\ref{fig:DRERM_Bench_gamma1}\subref{fig:DRERM_Bench_gamma1:gamma_2=1_others} indicate that the optimal solutions $x^*_{\gamma_1,1}$ to problem \eqref{prob:equivalent_NSDP2} are not changed for $1\leq\gamma_1\leq 2$.

To summarize the first experiment, as $\gamma_1$ increases, the solution $x^*_{\gamma_1,\gamma_2}$ tends to focus on decreasing the mean of realizations of $f_\alpha$ for the case of $\gamma_2=2$.
Moreover, the mean increases as $\gamma_1$ becomes larger when $\gamma_2=1$.
This implies that the uncertainty of the estimated variance-covariance $\tilde{\Sigma}_0$ is not sufficiently considered for the case of $\gamma_2=1$.

In the second experiment, we investigate the characteristics of the solutions in the case where $\gamma_2$ is incremented by 0.1 from 1 to 3, and $\gamma_1$ is set to 0.1 or 1.
We prepare 5000 realizations $\{\bar{\xi}^j\}_{j=1}^N$, which are the same samples used in the first experiment and compute the maximum, mean, and SD of $\{f_\alpha(x^*_{\gamma_1,\gamma_2},\bar{\xi}^j)\}_{j=1}^{N}$ for the solution $x^*_{\gamma_1,\gamma_2}$ to problem \eqref{prob:equivalent_NSDP2}.

\begin{figure}[htbt]
    \centering
    \subfloat[Maximum for fixed $\gamma_1=0.1$]{%
        \resizebox*{7cm}{!}{\includegraphics{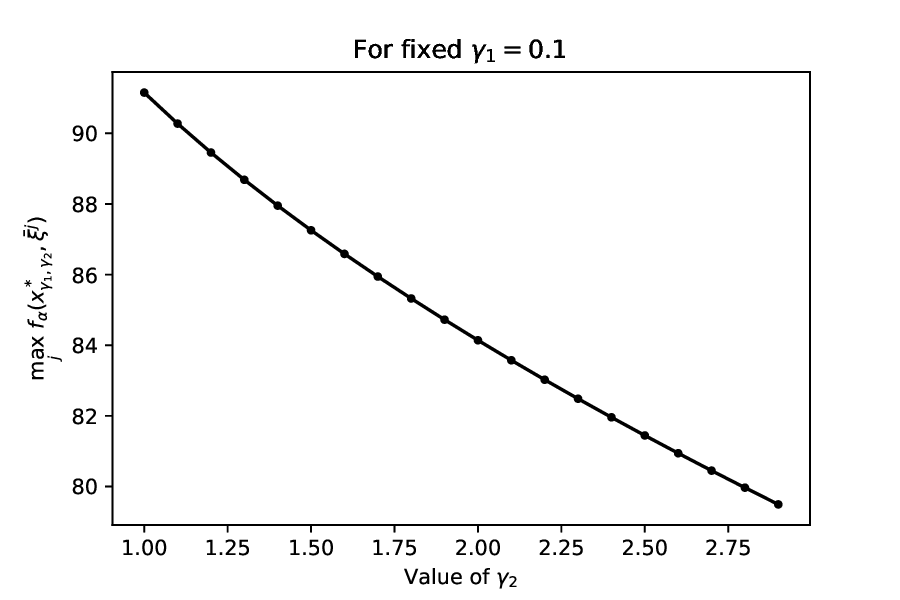}}\label{fig:DRERM_Bench_gamma2:gamma_1=0.1_max}
    }
    \subfloat[Maximum for fixed $\gamma_1=1$]{%
        \resizebox*{7cm}{!}{\includegraphics{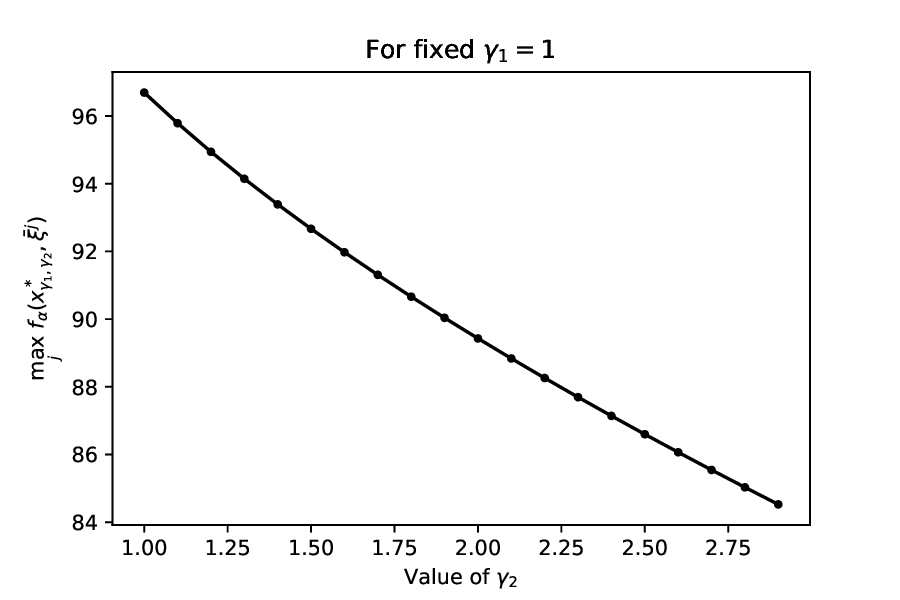}}\label{fig:DRERM_Bench_gamma2:gamma_1=1.0_max}
    }\\
    \subfloat[Mean and SD for fixed $\gamma_1=0.1$]{%
        \resizebox*{7cm}{!}{\includegraphics{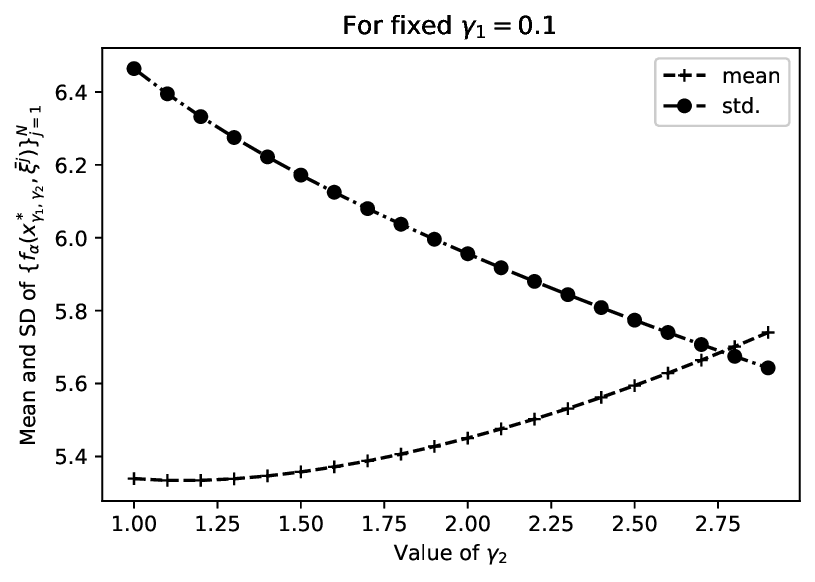}}\label{fig:DRERM_Bench_gamma2:gamma_1=0.1_others}
    }
    \subfloat[Mean and SD for fixed $\gamma_1=1$]{%
        \resizebox*{7cm}{!}{\includegraphics{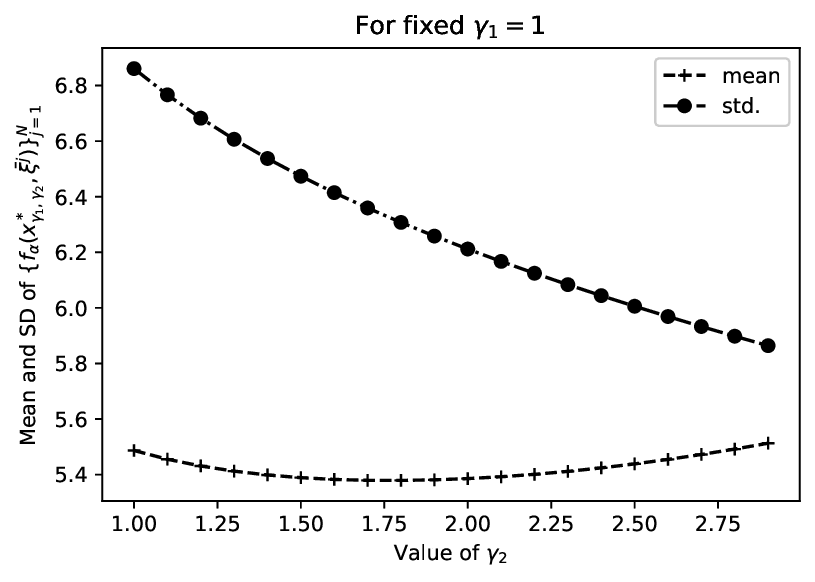}}\label{fig:DRERM_Bench_gamma2:gamma_1=1.0_others}
    }
    \caption{Maximum, mean, and SD of 5000 realizations of the regularized gap function when $\gamma_2$ is varied.}\label{fig:DRERM_Bench_gamma2}
\end{figure}

Figure~\ref{fig:DRERM_Bench_gamma2} depicts the results of the second experiment.
In particular, Figures~\ref{fig:DRERM_Bench_gamma2}\subref{fig:DRERM_Bench_gamma2:gamma_1=0.1_max} and~\ref{fig:DRERM_Bench_gamma2}\subref{fig:DRERM_Bench_gamma2:gamma_1=1.0_max} are $\max_j\ f_\alpha(x^*_{\gamma_1,\gamma_2},\bar{\xi}^j)$ for fixed $\gamma_1=0.1$ and $\gamma_1=1$, respectively.
Figures~\ref{fig:DRERM_Bench_gamma2}\subref{fig:DRERM_Bench_gamma2:gamma_1=0.1_others} and~\ref{fig:DRERM_Bench_gamma2}\subref{fig:DRERM_Bench_gamma2:gamma_1=1.0_others} are the mean and SD of $\{f_\alpha(x^*_{\gamma_1,\gamma_2},\bar{\xi}^j)\}_{j=1}^N$ for fixed $\gamma_1=0.1$ and $\gamma_1=1$, respectively.

For fixed $\gamma_1=0.1$, the maximum and SD gradually decrease as $\gamma_2$ increases, whereas the mean increases.
For fixed $\gamma_1=1$, the maximum and SD also decrease; however, the values of $f_\alpha$ are larger than the case where $\gamma_1=0.1$ entirely.
Moreover, there is diminutive change in the curve of the mean in Figure~\ref{fig:DRERM_Bench_gamma2}\subref{fig:DRERM_Bench_gamma2:gamma_1=1.0_others} compared with that of Figure~\ref{fig:DRERM_Bench_gamma2}\subref{fig:DRERM_Bench_gamma2:gamma_1=0.1_others}.

To summarize the second experiment, as $\gamma_2$ increases, the DRERM outputs the solutions $x^*_{\gamma_1,\gamma_2}$ that tend to decrease the maximum and SD of $f_\alpha$.
This is because, by the definition of the moment ambiguity set \eqref{eq:moment_set}, increasing $\gamma_2$ leads to the conservative behavior regarding the variance of $\xi$.
Consequently, $f_\alpha$ also behaves conservatively, and its outlier tends to be decreased as well.
Meanwhile, when $\gamma_2$ is very large, the mean increases.

Consequently, from the results of both the experiments, we confirm that there are trade-off relations between the mean and the SD, and the mean and the maximum, respectively, in response to the confidence parameters $\gamma_1$ and $\gamma_2$.

\begin{remark}
When the support $\Xi$ is compact, the reasonable $\gamma_1$ and $\gamma_2$ can be analytically obtained depending on the number of observations (refer to {\rm \cite{Delage2010}}).
However, if $\Xi$ is not compact, such as this experiment, one can obtain desired $\gamma_1$, $\gamma_2$, and solutions to SVIP \eqref{SVIP} by approximating $\Xi$ into a compact set.
\end{remark}

\section{Concluding remarks}\label{sec:conclusion}

We have proposed a DRERM model for an SVIP under uncertainty of distribution by incorporating the idea of the DRO into the ERM model with the regularized gap function.
In particular, we have shown that the DRERM can be conservatively approximated into a deterministic NSDP, and under suitable assumptions, the solution of the NSDP also solves the DRERM.
Furthermore, for the SVIP whose mapping $F$ is affine with respect to $x$, we have provided a sufficient condition of the regularization parameter of the regularized gap function to ensure that the reformulated NSDP is a convex optimization problem.
Meanwhile, the reformulated NSDP proposed in the existing research is not convex in general.
In numerical experiments, we have confirmed the reasonability of the DRERM model by comparing it with the ERM in terms of robustness, and we have analyzed their solutions by varying confidence parameters $\gamma_1$ and $\gamma_2$ included in the ambiguity set $\mathscr{P}$.


A remaining challenge is an NSDP approximation for more general cases of the following ambiguity sets described in \cite{Xu2018}:
\[
    \mathscr{P}'=\Curly{
        P\in\mathscr{M}_\Xi\ \middle|\ 
			\begin{array}{l}
                \ExpP{\Psi_i(\xi)}=O,\ i=1,2,\dots,t'\\
                \ExpP{\Psi_i(\xi)}\preceq O,\ i=t'+1,t'+2,\dots,t
			\end{array}
    },
\]
where $\Psi_i\ (i=1,2,\dots,t)$ is a symmetric matrix- or scalar-valued function over $\Xi$ with measurable random components.
We expect that our approach can be extended into the case of $\mathscr{P}'$ because the DRO with $\mathscr{P}'$ can be equivalently reformulated to a semi-infinite programming problem, such as (SIP), by assuming a `Slater-type' condition on $\mathscr{P}'$.

\section*{Acknowledgements}
The authors are grateful to two anonymous reviewers for careful reading of the manuscript and insightful comments to improve the quality of the paper.

\section*{Disclosure statement}
No potential conflict of interest was reported by the authors.

\section*{Funding}
This work was supported by the JSPS KAKENHI under Grant JP17K00032.


\appendix

\end{document}